\documentclass[11pt]{article}
\usepackage{}
\usepackage{amssymb}
\usepackage{amsfonts}
\usepackage{mathrsfs}
\usepackage{graphicx}
\usepackage{amsmath}
\usepackage{color}
\usepackage{amsthm}

\usepackage{pifont,bm}
\usepackage{tikz}
\usepackage{epstopdf}
\usepackage{enumerate}

\theoremstyle{definition}

\makeatletter
\def\@biblabel#1{[#1]}
\makeatother

\makeatletter \@addtoreset{equation}{section}

\begin{document}

\begin{titlepage}
\title{\bf{Long-time asymptotics for a fourth-order dispersive nonlinear Schr\"{o}dinger equation with nonzero boundary conditions
\footnote{
Corresponding author.\protect\\
\hspace*{3ex} \emph{E-mail addresses}: ychen@sei.ecnu.edu.cn (Y. Chen)}
}}
\author{Weiqi Peng$^{a}$, Yong Chen$^{a,b,*}$ \\
\small \emph{$^{a}$School of Mathematical Sciences, Shanghai Key Laboratory of PMMP} \\
\small \emph{East China Normal University, Shanghai, 200241, China} \\
\small \emph{$^{b}$College of Mathematics and Systems Science, Shandong University }\\
\small \emph{of Science and Technology, Qingdao, 266590, China} \\
\date{}}
\thispagestyle{empty}
\end{titlepage}
\maketitle

\vspace{-0.5cm}
\begin{center}
\rule{15cm}{1pt}\vspace{0.3cm}

\parbox{15cm}{\small
{\bf Abstract}\\
\hspace{0.5cm}  In this work, we consider the long-time asymptotics for the Cauchy problem of a fourth-order dispersive nonlinear Schr\"{o}dinger equation with nonzero boundary conditions at infinity. Firstly, in order to construct the basic Riemann-Hilbert  problem associated with nonzero boundary conditions,  we analysis  direct scattering problem. Then we deform  the corresponding matrix Riemann-Hilbert problem to explicitly solving models via using the  nonlinear steepest descent method and employing the $g$-function mechanism to eliminate the exponential growths of the jump matrices. Finally, we obtain the asymptotic stage of modulation instability  for the fourth-order dispersive nonlinear Schr\"{o}dinger equation.
}
\vspace{0.5cm}
\parbox{15cm}{\small{

\vspace{0.3cm} \emph{Key words:} Long-time asymptotics;  Fourth-order dispersive nonlinear Schr\"{o}dinger equation equation; Riemann-Hilbert problem; Nonlinear steepest descent method.\\

 } }
\end{center}
\vspace{0.3cm} \rule{15cm}{1pt} \vspace{0.2cm}

\tableofcontents

\vspace{0.3cm}

\section{Introduction}
The fourth-order dispersive nonlinear Schr\"{o}dinger(NLS) equation with nonzero boundary conditions(NZBCs) is
\begin{align}\label{1}
&iq_{t}+q_{xx}+2(|q|^{2}-q_{0}^{2})q+\gamma(q_{xxxx}+8q_{xx}|q|^{2}+2q_{xx}^{\ast}q^{2}\notag\\
&+6q_{x}^{2}q^{\ast}
+4q|q_{x}|^{2}+6(|q|^{4}-q_{0}^{4})q)=0,\notag\\
&\lim_{x\rightarrow\pm\infty}q(x, 0)=q_{\pm},
\end{align}
where $q_{\pm}$ are complex constants and independent of $x, t$ with $\mid q_{\pm}\mid=q_{0}>0$. Moreover
\begin{align}\label{1.1}
q(x, 0)-q_{\pm}\in L^{1,1}(\mathbb{R}^{\pm}),\ L^{1,1}(\mathbb{R}^{\pm})=\{f:\mathbb{R}\rightarrow \mathbb{C}\mid
\int_{\mathbb{R}^{\pm}}(1+|x|)|f(x)|dx<\infty\}.
\end{align}
With the Gauge transformation
\begin{align}\label{3}
q(x, t)=u(x, t)e^{-2i(3\gamma q_{0}^{2}+1)q_{0}^{2}t},
\end{align}
Eq.\eqref{1} can be reduced to the following  general fourth-order dispersive NLS equation\cite{Zhanghaiq}
\begin{align}\label{2}
&iu_{t}+u_{xx}+2|u|^{2}u+\gamma(u_{xxxx}+8u_{xx}|u|^{2}+2u_{xx}^{\ast}u^{2}+6u_{x}^{2}u^{\ast}
+4u|u_{x}|^{2}+6u|u|^{4})=0,\notag\\
&\lim_{x\rightarrow\pm\infty}u(x, t)=q_{\pm}e^{2i(3\gamma q_{0}^{2}+1)q_{0}^{2}t},
\end{align}
where $u$ is complex function with temporal variable $t$ and spatial variable $x$, which denotes the slowly varying envelope of the wave. The parameter $\gamma$ is a small dimensionless real numbers. In long distance and high speed optical fiber transmission systems, the fourth-order dispersion NLS equation plays a leading role in describing the transmission of ultrashort optical pulses \cite{Tian22,Tian23,Tian24}. Moreover, the equation can also depict the nonlinear spin excitation in a one-dimensional isotropic biquadratic Heisenberg ferromagnetic spin with octopole-dipole interaction\cite{Tian25,Tian26}. So far, there are some works on the study of the fourth-order dispersive NLS equation. Many methods have been used to derive the exact solutions for the fourth-order dispersive NLS equation, such as Darboux
transformation method, Hirota bilinear method and inverse scattering transform (IST) method\cite{Tian28,Tian-Li}, and the Lax pair, conservation laws, local wave solutions have also been
discussed\cite{guo-liu9,guo-liu18,guo-liu22}. Recently, the long-time asymptotic behavior for the fourth-order dispersive
NLS equation under zero boundary conditions(ZBCs) was investigated\cite{guo-liu,guo-liu11}. To our known of knowledge, the long-time asymptotic behaviors for the fourth-order dispersive NLS equation with NZBCs have not been analyzed yet.

In fact, the asymptotic behavior of solutions for nonlinear integrable systems has a long history and is always a hot topic. Early studies can be traced back to literatures\cite{Peng28,Peng29,Peng30,Peng31,Peng32,Peng33}. It is worth mentioning that Deift and Zhou, motivated by the pioneering work of Its\cite{Peng33}, proposed the nonlinear steepest descent method to investigate the long-time asymptotic behavior for the Cauchy problem of the mKdV equation with a oscillatory Riemann-Hilbert(RH) problem\cite{Peng34}. Subsequently,  this method was further developed in references\cite{Peng35,Peng36,Peng37}. Since the nonlinear steepest descent method been a efficient technique to research the Cauchy problem
of integrable equations,  the long-time asymptotics for lots of integrable equations as followed have been analyzed\cite{Peng38,Peng39,Peng40,Peng41,Peng43,Peng44,Wang-ke,Geng}. Besides, the method have been extended to the long-time asymptotics of the Cauchy problems for nonlinear integrable systems with a variety of non-decaying initial data, such as the time-periodic boundary conditions\cite{Peng46,Peng47}, the shock problem\cite{Peng45}, and  the step-like initial data\cite{Peng48,Peng49,Peng50}. Moreover, as a significant development of RH problem,   $\bar{\partial}$ generalization of the nonlinear steepest descent method was raised to derive the long-time asymptotic expansion of the solution in different fixed space-time regions\cite{Peng51,Peng52,Peng53,Tian-jde}.
Recent years, the researches about NZBCs at infinity  have
already been become a focal point. Biondini and his cooperators have studied the soliton solutions and the long-time asymptotics for the focusing NLS equation with NZBCs in \cite{Liunan5} and \cite{Liunan7}, respectively. After that, long-time asymptotics of the focusing Kundu-Eckhaus equation with NZBCs were studied in \cite{Wang}, long-time dynamics of the Gerdjikov-Ivanov type derivative nonlinear Schr\"{o}dinger equation with NZBCs were studied in \cite{Liunan}, long-time dynamics of the Hirota equation with NZBCs were studied in \cite{Chen-yan}, and long-time dynamics of the  modified Landau-Lifshitz equation with NZBCs were studied in \cite{Peng-tian}. Besides, the long-time asymptotic behavior of nonlocal integrable NLS solutions with NZBCs  were studied in \cite{Shepelsky}.

In this work, motivated by the long-time asymptotic analysis presented in \cite{Liunan7},  we consider the long-time asymptotics of Eq.\eqref{1} with the NZBCs at infinity. To the best knowledge of the authors, the long-time asymptotics for the fourth-order dispersive NLS equation under the NZBCs has never been reported up to now.

\textbf{The major results of this work is summarized in what follows:}

\noindent \textbf{Theorem 1.1.} \emph{
As $t\rightarrow\infty$, the asymptotic stage of modulation instability for $q(x,t)$ is given by
\begin{gather}
q(x,t)=\frac{q_{0}}{q_{+}^{\ast}}(\chi_{2}+q_{0})\frac{\Theta\left(\frac{\Omega t+\vartheta+i\ln\left(\frac{iq_{+}^{\ast}}{q_{0}}\right)}{2\pi}-V(\infty)+C\right)
\Theta(V(\infty)+C)}{\Theta\left(\frac{\Omega t+\vartheta+i\ln\left(\frac{iq_{+}^{\ast}}{q_{0}}\right)}{2\pi}+V(\infty)+C\right)\Theta(-V(\infty)+C)}
e^{2i(g(\infty)+G(\infty)t)}+\mathcal{O}(t^{-\frac{1}{2}}),
\end{gather}
where $V(\infty)=\int_{iq_{0}}^{\infty}d\vartheta$, and $\chi_{2}, \Omega, G(\infty), \vartheta, g(\infty), C$ are given by Eqs. \eqref{84}, \eqref{96}, \eqref{99}, \eqref{103a}, \eqref{103}, \eqref{117a}.
}

\textbf{Organization of this work:}  In Section 2, we perform the spectral analysis for the Cauchy problem of the fourth-order dispersive NLS equation, and construct the basic RH problem, which is the premise to give out the asymptotic behavior of the fourth-order dispersive NLS equation under NZBCs. In Section 3,  the asymptotic stage of modulation instability for Eq.\eqref{1} is analysed in detail.

\section{Reconstructing the basic Riemann-Hilbert problem}
We recall some notations that will used in our paper. The classical Pauli matrices are
defined as follows
\begin{align}\label{3.1}
\sigma_{1}=\left(\begin{array}{cc}
    0  &  1\\
    1 &  0\\
\end{array}\right),\qquad \sigma_{2}=\left(\begin{array}{cc}
    0  &  -i\\
    i &  0\\
\end{array}\right),\qquad \sigma_{3}=\left(\begin{array}{cc}
    1  &  0\\
    0 &  -1\\
\end{array}\right).
\end{align}
For a $2\times 2$ matrix $A$ and a scalar variable $\alpha$, we define
\begin{align}\label{3.2}
e^{\alpha\hat{\sigma}_{3}}A=e^{\alpha\sigma_{3}}Ae^{-\alpha\sigma_{3}}.
\end{align}

\subsection{Direct scattering problem with NZBCs}
The Lax pair of Eq.\eqref{1} is
\begin{align}\label{4}
\psi_{x}=X\psi,\qquad \qquad\qquad \psi_{t}=T\psi,
\end{align}
with the vector eigenfunction $\psi=(\psi_{1},\psi_{2})^{T}$ being a $2\times 2$ matrix, where the superscript $T$ represents the transpose of the vector,
\begin{align}\label{5}
X=&-ik\sigma_{3}+Q,\qquad
Q=\left(\begin{array}{cc}
    0  &  q\\
    -q^{\ast} &  0\\
\end{array}\right),\notag\\
T=&\left[3i\gamma|q|^{4}+i|q|^{2}+i\gamma(q_{xx}q^{\ast}+qq_{xx}^{\ast}-|q_{x}|^{2})-
2k\gamma(qq_{x}^{\ast}-q_{x}q^{\ast})\right.\notag\\
&\left.-2ik^{2}(2\gamma|q|^{2}+1)+8i\gamma k^{4}-i(3\gamma q_{0}^{2}+1)q_{0}^{2}\right]\sigma_{3}-4i\gamma k^{2}\sigma_{3}Q_{x}\notag\\
&-8\gamma k^{3}Q+6i\gamma Q^{2}Q_{x}\sigma_{3}+i\sigma_{3}Q_{x}+i\gamma\sigma_{3}Q_{xxx}+
2k(Q+\gamma Q_{xx}-2\gamma Q^{3}),
\end{align}
where $k$ represents the spectrum parameter. For convenience, we set $\gamma=1$ for the following analysis.

Taking $x\rightarrow\pm\infty$ and combining the NZBCs, we turn the Lax pair in Eq.\eqref{4} into
\begin{align}\label{6}
\psi_{\pm x}=X_{\pm}\psi_{\pm},\qquad \psi_{\pm t}=T_{\pm}\psi_{\pm},
\end{align}
where
\begin{gather}
X_{\pm}=-ik\sigma_{3}+Q_{\pm},\quad T_{\pm}=(-8k^{3}+2k+4kq_{0}^{2})X_{\pm}(k), \notag\\ Q_{\pm}=\left(\begin{array}{cc}
    0  &  q_{\pm}\\
    -q^{\ast}_{\pm} &  0\\
\end{array}\right),
\label{7}
\end{gather}
and we have defined  $Q_{\pm}=\lim\limits_{x\rightarrow \pm\infty}Q$.

It is not hard to obtain the eigenvalues of the matrix $X_{\pm}$ given by $\pm i\lambda$, and $\lambda=\sqrt{k^{2}+q_{0}^{2}}$. Obviously,
the branch cut of $\lambda$
is $\eta=[-iq_{0},iq_{0}]$ which is oriented upward in Figure 1. Here, we should also define  $\eta_{+}=[0,iq_{0}]$ and $\eta_{-}=[-iq_{0},0]$.\\

\centerline{\begin{tikzpicture}[scale=1.2]
\draw[->][thick](-3,0)--(-2,0);
\draw[-][thick](-2.0,0)--(0,0);
\draw[fill] (0.2,-0.2) node{$0$};
\draw[->][thick](0,0)--(2,0);
\draw[-][thick](2,0)--(3,0);
\draw[fill] (0,0) circle [radius=0.035];
\draw [->,very thick] (0,0)--(0,1);
\draw [-,very thick] (0,1)--(0,2);
\draw [->,very thick] (0,-2)--(0,-1);
\draw [-,very thick] (0,-1)--(0,0);
\draw[fill] (3,0.2) node{$\mathbb{R}$};
\draw[fill] (0.5,2) node{$iq_{0}$};
\draw[fill] (0.5,-2) node{$-iq_{0}$};
\draw[fill] (0.3,1) node{$\eta_{+}$};
\draw[fill] (0.3,-1) node{$\eta_{-}$};
\end{tikzpicture}}
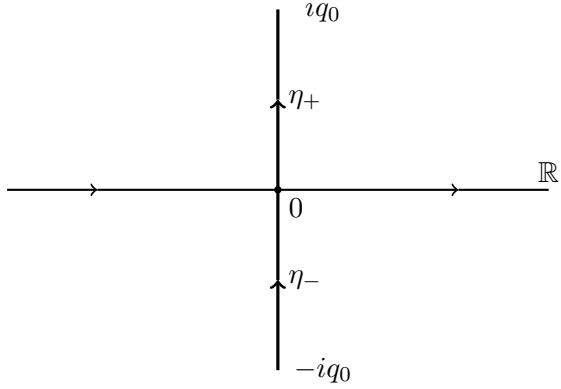
\centerline{\noindent {\small \textbf{Figure 1.} (Color online) The contour $\Sigma=\mathbb{R}\cup\eta$ of the basic RH problem.}}

The asymptotic spectral problem \eqref{6} can be solved by
\begin{align}\label{9}
\psi_{\pm}=E_{\pm}e^{-i\theta(x,t,k)\sigma_{3}},
\end{align}
where
\begin{align}\label{10}
\theta=\lambda\left[x+(-8k^{3}+2k+4kq_{0}^{2})t\right],\quad E_{\pm}=\left(\begin{array}{cc}
    1  &  \frac{\lambda-k}{iq^{\ast}_{\pm}}\\
    \frac{\lambda-k}{iq_{\pm}} &  1\\
\end{array}\right).
\end{align}

Supposing that $\Psi_{\pm}(x,t,k)$ are both the Jost solutions of the Lax pair \eqref{4}, we can define $\Psi_{\pm}(x,t,k)=\psi_{\pm}(x,t,k)+o(1)$ as $x\rightarrow\infty$. Furthermore, using variable transformation
\begin{align}\label{11}
\mu_{\pm}(x,t,k)=\Psi_{\pm}(x,t,k)e^{i\theta(x,t,k)\sigma_{3}},
\end{align}
we have
\begin{align}\label{12}
\mu_{\pm}(x,t,k)=E_{\pm}+o(1),\qquad x\rightarrow\pm\infty,
\end{align}
which arrive the following two Volterra integral equations
\begin{align}\label{13}
\mu_{-}(x,t,k)=E_{-}+\int_{-\infty}^{x}E_{-}e^{-i\lambda(x-y)\hat{\sigma}_{3}}
\left[E^{-1}_{-}(Q-Q_{-})\mu_{-}(y,t,k)\right]dy,\notag\\
\mu_{+}(x,t,k)=E_{+}-\int_{x}^{+\infty}E_{+}e^{-i\lambda(x-y)\hat{\sigma}_{3}}
\left[E^{-1}_{+}(Q-Q_{+})\mu_{+}(y,t,k)\right]dy.
\end{align}

\noindent \textbf{Proposition 2.1.} \emph{
Suppose $q-q_{\pm}\in L^{1} (\mathbb{R}^{\pm})$, then $\mu_{\pm}(x,t,k)$ given in Eq.\eqref{12} uniquely satisfy the Volterra integral equation \eqref{13} in $\Sigma_{0}$, and $\mu_{\pm}(x,t,k)$ admit:}

\emph{$\bullet$ $\mu_{-1}(x, t, k)$ and $\mu_{+2}(x, t, k)$ is analytical in $\mathbb{C}_{+}\setminus \eta_{+}$ and
continuous in $\mathbb{C}_{+}\cup \Sigma_{0}$;}

\emph{$\bullet$ $\mu_{+1}(x, t, k)$ and $\mu_{-2}(x, t, k)$  is analytical in $\mathbb{C}_{-}\setminus \eta_{-}$ and
continuous in $\mathbb{C}_{-}\cup \Sigma_{0}$;}

\emph{$\bullet$ $\mu_{\pm}(x,t,k)\rightarrow I$ \mbox{as}  $k\rightarrow \infty$;}

\emph{$\bullet$ $\det \mu_{\pm}(x, t, k)=\det E_{\pm}=\frac{2\lambda}{\lambda+ k}\triangleq d(k), \quad x, t\in \mathbb{R}, \quad k\in \Sigma_{0}$.}

Since Jost solutions $\Psi_{\pm}(x,t,k)$ meet the Lax pair \eqref{4} for $k\in \Sigma_{0}$, the linear relation can be established by scattering matrix $s(k)$
\begin{align}\label{15}
\Psi_{-}(x,t,k)=\Psi_{+}(x,t,k)s(k),\qquad k\in \Sigma_{0},
\end{align}
of which the scattering matrix $s(k)$ has the following symmetry properties
\begin{align}\label{16}
\Psi^{\ast}_{\pm}(k^{\ast})=\sigma_{2}\Psi_{\pm}(k)\sigma_{2},\ s^{\ast}(k^{\ast})=\sigma_{2}s(k)\sigma_{2}, \ k\in \Sigma_{0}.
\end{align}
Therefore, the concrete form of scattering matrix $s(k)$ arrives at
\begin{align}\label{17}
s(k)=\left(\begin{array}{cc}
   a(k)  &  -b^{\ast}(k)\\
    b(k) &  a^{\ast}(k)\\
\end{array}\right),\qquad a(k)a^{\ast}(k)+b(k)b^{\ast}(k)=1,
\end{align}
where $a^{\ast}(k)=a^{\ast}(k^{\ast}), b^{\ast}(k)=b^{\ast}(k^{\ast})$ means the Schwartz conjugates. Then we get
\begin{align}\label{18}
a(k)=\frac{Wr(\Psi_{-1},\Psi_{+2})}{d(k)}, \quad a^{\ast}(k)=\frac{Wr(\Psi_{+1},\Psi_{-2})}{d(k)},\notag\\
b(k)=\frac{Wr(\Psi_{+1},\Psi_{-1})}{d(k)}, \quad b^{\ast}(k)=\frac{Wr(\Psi_{+2},\Psi_{-2})}{d(k)}.
\end{align}
By taking $\eta$ to be oriented upwards and defining
\begin{align}\label{19}
\begin{aligned}
&\mu_{-1}^{-}(k)=\lim\limits_{\varepsilon\rightarrow 0^{+}}\mu_{-1}(k+\varepsilon)=\mu_{-1}(k),\\
&\mu_{+2}^{-}(k)=\lim\limits_{\varepsilon\rightarrow 0^{+}}\mu_{+2}(k+\varepsilon)=\mu_{+2}(k),
\end{aligned}
\qquad\qquad k\in \eta_{-}.
\end{align}
we derive  the Jost
solutions $\mu_{\pm}$ and the scattering data $a, b$ have the following jump conditions across the
branch cut $\eta$ respectively, given by
\begin{align}\label{20}
\begin{aligned}
&\mu_{+1}^{+}(x,t,k)=\frac{i(\lambda+ k)}{q_{+}}\mu_{+2}(x,t,k),\\
&\mu_{-2}^{+}(x,t,k)=\frac{i(\lambda+ k)}{q^{\ast}_{-}}\mu_{-1}(x,t,k),
\end{aligned}
\qquad k\in \eta_{+}. \notag\\
\begin{aligned}
&\mu_{-1}^{+}(x,t,k)=\frac{i(\lambda+ k)}{q_{-}}\mu_{-2}(x,t,k),\\
&\mu_{+2}^{+}(x,t,k)=\frac{i(\lambda+ k)}{q^{\ast}_{+}}\mu_{+1}(x,t,k),
\end{aligned}
\qquad k\in \eta_{-},
\end{align}
and
\begin{align}\label{21}
(a^{\ast})^{+}(k)=\frac{q_{-}}{q_{+}}a(k),\ a^{+}(k)=\frac{q_{+}}{q_{-}}a^{\ast}(k),\ b^{+}(k)=-\frac{q_{-}^{\ast}}{q_{+}}b^{\ast}(k).
\end{align}
\subsection{Inverse scattering problem and reconstructing the formula for potential}
Now, the fundamental matrix-value function is formulated as
\begin{align}\label{22}
m(x, t, k)=\left\{
\begin{aligned}
&(\frac{\Psi_{-1}}{ad},\Psi_{+2})e^{i\theta\sigma_{3}}, \qquad k\in\mathbb{C}_{+}\setminus \eta_{+},\\
&(\Psi_{+1},\frac{\Psi_{-2}}{a^{\ast}d})e^{i\theta\sigma_{3}}, \qquad k\in\mathbb{C}_{-}\setminus \eta_{-},
  \end{aligned}
\right.
\end{align}
Then the matrix-value function $m(x, t, k)$ has following jump condition across $\mathbb{R}$:
\begin{align}\label{23}
m_{+}(x, t, k)=m_{-}(x, t, k)\left(\begin{array}{cc}
   \frac{1}{d}[1+\gamma(k)\gamma^{\ast}(k)]  &  \gamma^{\ast}(k)e^{-2i\theta(x,t,k)}\\
    \gamma(k)e^{2i\theta(x,t,k)} &  d(k)\\
\end{array}\right),\qquad k\in \mathbb{R},
\end{align}
where $m_{\pm}(x, t, k)$ denote the boundary values of $m(x, t, k)$ as $k$ approaches the contour
from  a chosen side, and the reflection coefficient $\gamma(k)=\frac{b(k)}{a(k)}$.
In terms of \eqref{21}, \eqref{22}, \eqref{23},
the jump condition of the matrix-value function $m(x, t, k)$ across $\eta_{+}$ is given by
\begin{align}\label{24}
m_{+}(x, t, k)=m_{-}(x, t, k)\left(\begin{array}{cc}
   \frac{\lambda-k}{iq_{+}}\gamma^{\ast}(k)e^{-2i\theta(x,t,k)}  &  \frac{2i\lambda}{q^{\ast}_{+}}\\
     \frac{iq^{\ast}_{+}}{2\lambda}[1+\gamma(k)\gamma^{\ast}(k)]&   \frac{\lambda+k}{iq^{\ast}_{+}}\gamma(k)e^{2i\theta(x,t,k)}\\
\end{array}\right),\ k\in \eta_{+}.
\end{align}
Similarly, the  matrix-value function $m^{(0)}(x, t, k)$ has jump condition across $\eta_{-}$:
\begin{align}\label{25}
m_{+}(x, t, k)=m_{-}(x, t, k)\left(\begin{array}{cc}
   \frac{i(\lambda+k)}{q_{+}}\gamma^{\ast}(k)e^{-2i\theta(x,t,k)}  &   \frac{iq_{+}}{2\lambda}[1+\gamma(k)\gamma^{\ast}(k)]\\
  \frac{2i\lambda}{iq_{+}}  &  \frac{i(\lambda-k)}{q^{\ast}_{+}}\gamma(k)e^{2i\theta(x,t,k)}\\
\end{array}\right),\ k\in \eta_{-}.
\end{align}
Finally, assuming that the $a\neq 0$ for all $k\in \mathbb{C}_{+}\cup \Sigma$, then a matrix RH
problem is constructed:
\begin{align}\label{26}
\left\{
\begin{array}{lr}
m(x, t, k)\ \mbox{is analytic in} \ \mathbb{C}\setminus\Sigma,\\
m_{+}(x, t, k)=m_{-}(x, t, k)J(x, t, k), \qquad k\in\Sigma,\\
m(x, t, k)\rightarrow I,\qquad k\rightarrow \infty,
  \end{array}
\right.
\end{align}
of which the jump matrix $J(x, t, k)=\{J_{i}(x, t, k)\}_{i=1}^{3}$ is (see Figure 1)
\begin{gather}
J_{1}=\left(\begin{array}{cc}
   \frac{1}{d}[1+\gamma(k)\gamma^{\ast}(k)]  &  \gamma^{\ast}(k)e^{-2if(x,t,k)t}\\
    \gamma(k)e^{2if(x,t,k)t} &  d(k)\\
\end{array}\right),\notag\\
J_{2}=\left(\begin{array}{cc}
   \frac{\lambda-k}{iq_{+}}\gamma^{\ast}(k)e^{-2if(x,t,k)t}  &  \frac{2i\lambda}{q^{\ast}_{+}}\\
     \frac{iq^{\ast}_{+}}{2\lambda}[1+\gamma(k)\gamma^{\ast}(k)]&   \frac{\lambda+k}{iq^{\ast}_{+}}\gamma(k)e^{2if(x,t,k)t}\\
\end{array}\right),\notag\\
J_{3}=\left(\begin{array}{cc}
   \frac{i(\lambda+k)}{q_{+}}\gamma^{\ast}(k)e^{-2if(x,t,k)t}  &   \frac{iq_{+}}{2\lambda}[1+\gamma(k)\gamma^{\ast}(k)]\\
  \frac{2i\lambda}{iq_{+}}  &  \frac{i(\lambda-k)}{q^{\ast}_{+}}\gamma(k)e^{2if(x,t,k)t}\\
\end{array}\right),\notag
\end{gather}
where $f=\lambda\left[\xi-8k^{3}+2k+4kq_{0}^{2}\right], \xi=\frac{x}{t}$.

In addition, expanding the $M^{(0)}(x, t, k)$ at large $k$ as
\begin{align}\label{27}
m(x, t, k)=I+\frac{m_{1}(x, t)}{k}+\frac{m_{2}(x, t)}{k^{2}}+\mathcal{O}(\frac{1}{k^{3}}),\qquad k\rightarrow \infty,
\end{align}
and combining equations \eqref{4}, \eqref{22}, \eqref{27}, we recover the solution $q(x, t)$ of the original initial value problem
\eqref{1} in the following form
\begin{align}\label{28}
q(x,t)=2i\left(m_{1}(x,t)\right)_{12}=2i\lim_{k\rightarrow\infty}km_{12}(x,t,k).
\end{align}

\subsection{The sign structure of $\mbox{Re}(if)$}
To find contour deformations, we firstly need to discuss the sign structure of the quantity $\mbox{Re}(if)$. Through taking
\begin{align}\label{36}
\frac{d f(\xi,k)}{dk}=\frac{-32k^{4}-(16q_{0}^{2}-4)k^{2}+k\xi+4q_{0}^{4}+2q_{0}^{2}}{\sqrt{k^{2}+q_{0}^{2}}},
\end{align}
for convenience, let $q_{0}^{2}=\frac{1}{2}$, we can get the four stationary phase points(i.e., the points $k_{s}$ such that $f'_{k}(k_{s})=0$), given by
\begin{align}\label{37}
k_{1}=\frac{-\sqrt{2y}\pm\sqrt{-\frac{\xi}{16\sqrt{2y}}-2y}}{2}, \quad k_{2}=\frac{\sqrt{2y}\pm\sqrt{\frac{\xi}{16\sqrt{2y}}-2y}}{2},
\end{align}
where
\begin{align}\label{38}
y=\left(\frac{\sqrt{\xi^{4}+128}}{16384}+\frac{\xi^{2}}{16384}\right)^{\frac{1}{3}}
-\left(\frac{\sqrt{\xi^{4}+128}}{16384}-\frac{\xi^{2}}{16384}\right)^{\frac{1}{3}}.
\end{align}
Using Maple symbolic computation, we find there are two real stationary phase points and two complex stationary phase points for arbitrary $\xi\neq 0$, and the corresponding sign structure of $\mbox{Re}(if)$ is presented at what follows. \\
\centerline{\begin{tikzpicture}[scale=1.5]
\draw[-][thick](-3,0)--(3,0);
\draw[-][thick](0,-1)--(0,1);
\draw [-,thick] (-2,0) to [out=90,in=-60] (-2.9,1);
\draw [-,thick] (-2,0) to [out=-90,in=60] (-2.9,-1);
\draw [-,thick] (2,0) to [out=90,in=-120] (2.9,1);
\draw [-,thick] (2,0) to [out=-90,in=120] (2.9,-1);
\draw [-,thick] (0,1) to [out=180,in=-60] (-0.5,2);
\draw [-,thick] (0,-1) to [out=180,in=60] (-0.5,-2);
\draw[fill] (-2.8,0.25) node{$\mathfrak{R}(if)<0$};
\draw[fill] (-2.8,-0.25) node{$\mathfrak{R}(if)>0$};
\draw[fill] (-1.5,1) node{$\mathfrak{R}(if)>0$};
\draw[fill] (-1.5,-1) node{$\mathfrak{R}(if)<0$};
\draw[fill] (1.5,1) node{$\mathfrak{R}(if)<0$};
\draw[fill] (1.5,-1) node{$\mathfrak{R}(if)>0$};
\draw[fill] (2.8,0.25) node{$\mathfrak{R}(if)>0$};
\draw[fill] (2.8,-0.25) node{$\mathfrak{R}(if)<0$};
\end{tikzpicture}}\\
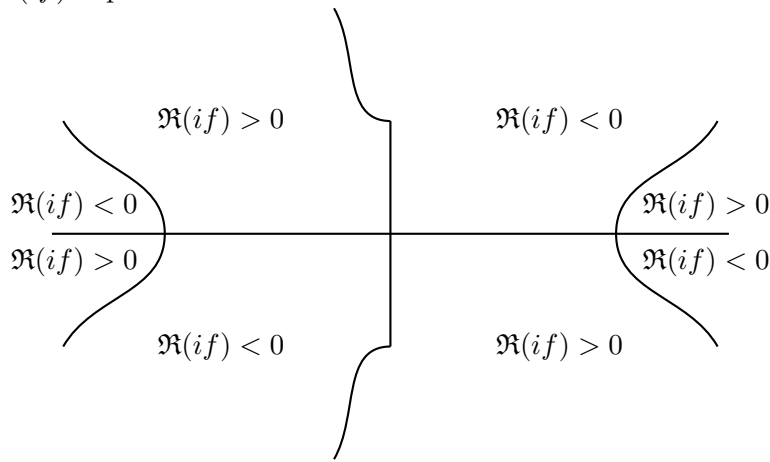
\noindent { \small \textbf{Figure 2.} (Color online) The sign structure of $\mbox{Re}(if)$ in the complex $k$-plane.}\\

\section{Asymptotic stage of modulation instability}
In order to derive long-time asymptotics of solution for the Eq.\eqref{1}, we carry out similar deformations of the RH problem \eqref{26} as that in Refs. \cite{Peng45,Liunan7}. As shown in Figure 2, the curves $\mbox{Im}\theta(k)=0$ will not intersect the real axis.  In order to study the long-time asymptotics of $q(x,t)$ in this case, we construct a $g$-function mechanism and introduce the point $k_{0}$.

\subsection{First deformation}
To achieve the first deformation, we first decompose the jump matrix $J_{1}, J_{2}, J_{3}$ into following form
\begin{align}
\left\{
\begin{array}{lr}
J_{1}=J_{2}^{(1)}J_{0}^{(1)}J_{1}^{(1)},\qquad \quad \ \mbox{on} \qquad (k_{1}, k_{0})\cup(k_{2}, \infty),\\
J_{1}=J_{4}^{(1)}J_{3}^{(1)},\qquad\qquad \quad \mbox{on} \qquad (-\infty, k_{1})\cup(k_{0}, k_{2}),\\
J_{2}=(J_{3-}^{(1)})^{-1}J_{\eta}^{(1)}J_{3+}^{(1)},\qquad \mbox{on} \qquad \eta_{+} \qquad \mbox{cut},\\
J_{3}=J_{4-}^{(1)}J_{\eta}^{(1)}(J_{4+}^{(1)})^{-1},\qquad \mbox{on} \qquad \eta_{-} \qquad \mbox{cut},\notag
\end{array}
\right.
\end{align}
where
\begin{gather}
J_{0}^{(1)}=\left(\begin{array}{cc}
   1+\gamma\gamma^{\ast}  &  0\\
   0 &  \frac{1}{1+\gamma\gamma^{\ast}}\\
\end{array}\right),\
J_{1}^{(1)}=\left(\begin{array}{cc}
   d^{-\frac{1}{2}}  &  \frac{d^{\frac{1}{2}}\gamma^{\ast}e^{-2if t}}{1+\gamma\gamma^{\ast}}\\
    0&  d^{\frac{1}{2}}\\
\end{array}\right),\ J_{2}^{(1)}=\left(\begin{array}{cc}
   d^{-\frac{1}{2}}  &  0\\
    \frac{d^{\frac{1}{2}}\gamma e^{2if t}}{1+\gamma\gamma^{\ast}}&  d^{\frac{1}{2}}\\
\end{array}\right),
\notag\\
J_{3}^{(1)}=\left(\begin{array}{cc}
   d^{-\frac{1}{2}}  &   0\\
  d^{-\frac{1}{2}}\gamma e^{2if t}  &  d^{\frac{1}{2}}\\
\end{array}\right),\ J_{4}^{(1)}=\left(\begin{array}{cc}
   d^{-\frac{1}{2}}  &   d^{-\frac{1}{2}}\gamma^{\ast} e^{-2if t}\\
    0 &  d^{\frac{1}{2}}\\
\end{array}\right),\ J_{\eta}^{(1)}=\left(\begin{array}{cc}
   0  &   \frac{iq_{+}}{q_{0}}\\
  \frac{iq^{\ast}_{+}}{q_{0}} &  0\\
\end{array}\right). \label{40}
\end{gather}
Therefore, we can transform $m$ into $m^{(1)}$ by using
\begin{align}\label{41}
m^{(1)}=mB(k),
\end{align}
where
\begin{align}\label{42}
B(k)=\left\{
\begin{array}{lr}
(J_{1}^{(1)})^{-1}\ \mbox{on} \ k\in\Omega_{1},\\
J_{2}^{(1)}\ \mbox{on} \ k\in\Omega_{2},\\
(J_{3}^{(1)})^{-1}\ \mbox{on} \ k\in\Omega_{3}\cup\Omega_{5},\\
J_{4}^{(1)}\ \mbox{on} \ k\in\Omega_{4}\cup\Omega_{6},\\
I \ \mbox{on} \ k\in \mbox{others},
  \end{array}
\right.
\end{align}
then the following RH problem about $m^{(1)}$ can be given
\begin{align}\label{42.1}
\left\{
\begin{array}{lr}
m^{(1)}(x, t, k)\ \mbox{is analytic in} \ \mathbb{C}\setminus\Sigma^{(1)},\\
m^{(1)}_{+}(x, t, k)=m^{(1)}_{-}(x, t, k)J^{(1)}(x, t, k), \qquad k\in\Sigma^{(1)},\\
m^{(1)}(x, t, k)\rightarrow I,\qquad k\rightarrow \infty,
  \end{array}
\right.
\end{align}
where the jump matrix $J^{(1)}$ has been defined in \eqref{40}.
The contour $\Sigma^{(0)}$  in Figure 3 also become
the  new contour $\Sigma^{(1)}$ as shown in Figure 4.\\

\centerline{\begin{tikzpicture}
\draw[-][dashed,thick](-6,0)--(-4,0);
\draw[-][dashed,thick](-2,0)--(3,0);
\draw[->][thick](-4,0)--(-3,0);
\draw[-][thick](-3,0)--(-2,0);
\draw[->][thick](3,0)--(4,0);
\draw[-][thick](4,0)--(6,0);
\draw[->][thick](0,0)--(0,1);
\draw[-][thick](0,1)--(0,2);
\draw[-][thick](0,0)--(0,-1);
\draw[<-][thick](0,-1)--(0,-2);
\draw [-,dashed,thick] (-4,0) to [out=90,in=-60] (-5,2);
\draw [-,dashed,thick] (-4,0) to [out=-90,in=60] (-5,-2);
\draw [-,dashed,thick] (3,0) to [out=90,in=-120] (4,2);
\draw [-,dashed,thick] (3,0) to [out=-90,in=120] (4,-2);
\draw [-,dashed,thick] (0,2) to [out=180,in=-60] (-0.5,3);
\draw [-,dashed,thick] (0,-2) to [out=180,in=60] (-0.5,-3);
\draw[->][thick](-6,1)--(-5,1);
\draw [-,thick] (-5,1) to [out=-5,in=145] (-4,0);
\draw[->][thick](-6,-1)--(-5,-1);
\draw [-,thick] (-5,-1) to [out=5,in=-145] (-4,0);
\draw [->,thick] (-4,0) to [out=35,in=180] (-3,1);
\draw [-,thick] (-3,1) to [out=0,in=145] (-2,0);
\draw [->,thick] (-4,0) to [out=-35,in=-180] (-3,-1);
\draw [-,thick] (-3,-1) to [out=0,in=-145] (-2,0);
\draw [->,thick] (-2,0) to [out=35,in=180] (-1,0.6);
\draw [-,thick] (-1,0.6) to [out=0,in=145] (0,0);
\draw [->,thick] (-2,0) to [out=-35,in=-180] (-1,-0.6);
\draw [-,thick] (-1,-0.6) to [out=0,in=-145] (0,0);
\draw [->,thick] (0,0) to [out=35,in=180] (1.5,1);
\draw [-,thick] (1.5,1) to [out=0,in=145] (3,0);
\draw [->,thick] (0,0) to [out=-35,in=-180] (1.5,-1);
\draw [-,thick] (1.5,-1) to [out=0,in=-145] (3,0);
\draw [->,thick] (3,0) to [out=35,in=180] (4,1);
\draw[-][thick](4,1)--(6,1);
\draw [->,thick] (3,0) to [out=-35,in=-180] (4,-1);
\draw[-][thick](4,-1)--(6,-1);
\draw [->,thick] (0,0) to [out=50,in=-90] (0.6,1);
\draw [-,thick]  (0.6,1) to [out=90,in=-45] (0,2);
\draw [->,thick] (0,0) to [out=130,in=-90] (-0.6,1);
\draw [-,thick]  (-0.6,1) to [out=90,in=-135] (0,2);
\draw [-,thick] (0,0) to [out=-50,in=90] (0.6,-1);
\draw [<-,thick]  (0.6,-1) to [out=-90,in=45] (0,-2);
\draw [-,thick] (0,0) to [out=-130,in=90] (-0.6,-1);
\draw [<-,thick]  (-0.6,-1) to [out=-90,in=135] (0,-2);
\draw[fill] (3,0) node[below]{$k_{2}$};
\draw[fill] (-2,0) node[below]{$k_{0}$};
\draw[fill] (-4,0) node[below]{$k_{1}$};
\draw[fill] (-6,1) node[above]{$\textcolor[rgb]{1.00,0.00,0.00}{J_{3}^{(1)}}$};
\draw[fill] (-6,-1) node[below]{$\textcolor[rgb]{1.00,0.00,0.00}{J_{4}^{(1)}}$};
\draw[fill] (-3,1) node[above]{$\textcolor[rgb]{1.00,0.00,0.00}{J_{1}^{(1)}}$};
\draw[fill] (-3,-1) node[below]{$\textcolor[rgb]{1.00,0.00,0.00}{J_{2}^{(1)}}$};
\draw[fill] (-2.5,0) node{$\textcolor[rgb]{1.00,0.00,0.00}{J_{0}^{(1)}}$};
\draw[fill] (-1.5,0.6) node{$\textcolor[rgb]{1.00,0.00,0.00}{J_{3}^{(1)}}$};
\draw[fill] (-1.5,-0.6) node{$\textcolor[rgb]{1.00,0.00,0.00}{J_{4}^{(1)}}$};
\draw[fill] (2.3,0.6) node[above]{$\textcolor[rgb]{1.00,0.00,0.00}{J_{3}^{(1)}}$};
\draw[fill] (2.3,-0.6) node[below]{$\textcolor[rgb]{1.00,0.00,0.00}{J_{4}^{(1)}}$};
\draw[fill] (5,1) node[above]{$\textcolor[rgb]{1.00,0.00,0.00}{J_{1}^{(1)}}$};
\draw[fill] (5,-1) node[below]{$\textcolor[rgb]{1.00,0.00,0.00}{J_{2}^{(1)}}$};
\draw[fill] (5,0) node{$\textcolor[rgb]{1.00,0.00,0.00}{J_{0}^{(1)}}$};
\draw[fill] (-0.4,1.3) node[left]{$\textcolor[rgb]{1.00,0.00,0.00}{J_{3}^{(1)}}$};
\draw[fill] (0.4,1.3) node[right]{$\textcolor[rgb]{1.00,0.00,0.00}{(J_{3}^{(1)})^{-1}}$};
\draw[fill] (-0.4,-1.3) node[left]{$\textcolor[rgb]{1.00,0.00,0.00}{(J_{4}^{(1)})^{-1}}$};
\draw[fill] (0.4,-1.3) node[right]{$\textcolor[rgb]{1.00,0.00,0.00}{J_{4}^{(1)}}$};
\draw[fill] (0,0.75) node{$\textcolor[rgb]{1.00,0.00,0.00}{J_{\eta}^{(1)}}$};
\draw[fill] (0,-0.75) node{$\textcolor[rgb]{1.00,0.00,0.00}{J_{\eta}^{(1)}}$};
\draw[fill] (-5,0.45) node[left]{$\Omega_{3}$};
\draw[fill] (-5,-0.45) node[left]{$\Omega_{4}$};
\draw[fill] (-3,0.45) node{$\Omega_{1}$};
\draw[fill] (-3,-0.45) node{$\Omega_{2}$};
\draw[fill] (-1,0.25) node{$\Omega_{3}$};
\draw[fill] (-1,-0.25) node{$\Omega_{4}$};
\draw[fill] (1.5,0.25) node{$\Omega_{3}$};
\draw[fill] (1.5,-0.25) node{$\Omega_{4}$};
\draw[fill] (5,0.45) node[left]{$\Omega_{1}$};
\draw[fill] (5,-0.45) node[left]{$\Omega_{2}$};
\draw[fill] (-0.3,1.25) node{$\Omega_{5}$};
\draw[fill] (-0.3,-1.25) node{$\Omega_{6}$};
\draw[fill] (0.3,1.25) node{$\Omega_{5}$};
\draw[fill] (0.3,-1.25) node{$\Omega_{6}$};
\end{tikzpicture}}
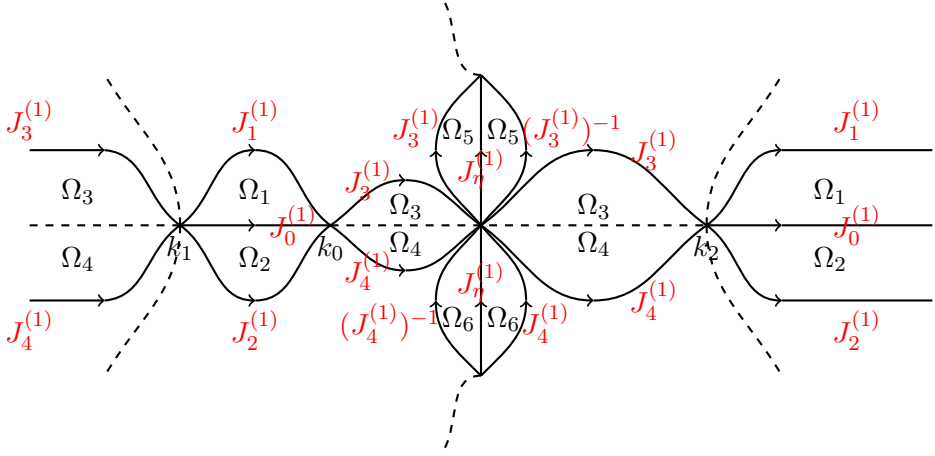
\centerline{\noindent {\small \textbf{Figure 3.} (Color online) The initial  contour  $\Sigma^{(0)}$.}}

\centerline{\begin{tikzpicture}
\draw[->][thick](-4,0)--(-3,0);
\draw[-][thick](-3,0)--(-2,0);
\draw[->][thick](3,0)--(4,0);
\draw[-][thick](4,0)--(6,0);
\draw[-][thick](0,0)--(0,1);
\draw[-][thick](0,1)--(0,1.4);
\draw[<-][thick](0,0)--(0,-1);
\draw[-][thick](0,-1)--(0,-1.4);
\draw [-,dashed,thick] (-4,0) to [out=90,in=-60] (-5,2);
\draw [-,dashed,thick] (-4,0) to [out=-90,in=60] (-5,-2);
\draw [-,dashed,thick] (3,0) to [out=90,in=-120] (4,2);
\draw [-,dashed,thick] (3,0) to [out=-90,in=120] (4,-2);
\draw [-,dashed,thick] (0,1.4) to [out=180,in=-60] (-0.5,3);
\draw [-,dashed,thick] (0,-1.4) to [out=180,in=60] (-0.5,-3);
\draw[->][thick](-6,1)--(-5,1);
\draw [-,thick] (-5,1) to [out=-5,in=145] (-4,0);
\draw[->][thick](-6,-1)--(-5,-1);
\draw [-,thick] (-5,-1) to [out=5,in=-145] (-4,0);
\draw [->,thick] (-4,0) to [out=35,in=180] (-3,1);
\draw [-,thick] (-3,1) to [out=0,in=145] (-2,0);
\draw [->,thick] (-4,0) to [out=-35,in=-180] (-3,-1);
\draw [-,thick] (-3,-1) to [out=0,in=-145] (-2,0);
\draw [->,thick] (3,0) to [out=35,in=180] (4,1);
\draw[-][thick](4,1)--(6,1);
\draw [->,thick] (3,0) to [out=-35,in=-180] (4,-1);
\draw[-][thick](4,-1)--(6,-1);
\draw[fill] (3,0) node[below]{$k_{2}$};
\draw[fill] (-2,0) node[below]{$k_{0}$};
\draw[fill] (-4,0) node[below]{$k_{1}$};
\draw[fill] (-6,1) node[above]{$\textcolor[rgb]{1.00,0.00,0.00}{J_{3}^{(1)}}$};
\draw[fill] (-6,-1) node[below]{$\textcolor[rgb]{1.00,0.00,0.00}{J_{4}^{(1)}}$};
\draw[fill] (-3,1) node[above]{$\textcolor[rgb]{1.00,0.00,0.00}{J_{1}^{(1)}}$};
\draw[fill] (-3,-1) node[below]{$\textcolor[rgb]{1.00,0.00,0.00}{J_{2}^{(1)}}$};
\draw[fill] (-2.6,0) node{$\textcolor[rgb]{1.00,0.00,0.00}{J_{0}^{(1)}}$};
\draw[fill] (5,1) node[above]{$\textcolor[rgb]{1.00,0.00,0.00}{J_{1}^{(1)}}$};
\draw[fill] (5,-1) node[below]{$\textcolor[rgb]{1.00,0.00,0.00}{J_{2}^{(1)}}$};
\draw[fill] (5,0) node{$\textcolor[rgb]{1.00,0.00,0.00}{J_{0}^{(1)}}$};
\draw[fill] (0.3,0.0) node{$\textcolor[rgb]{1.00,0.00,0.00}{J_{\eta}^{(1)}}$};
\draw [->,thick] (-2,0) to [out=35,in=180] (0.5,2);
\draw [-,thick] (0.5,2) to [out=0,in=145] (3,0);
\draw [->,thick] (-2,0) to [out=-35,in=-180] (0.5,-2);
\draw [-,thick] (0.5,-2) to [out=0,in=-145] (3,0);
\draw[fill] (-1,1) node[above]{$\textcolor[rgb]{1.00,0.00,0.00}{J_{3}^{(1)}}$};
\draw[fill] (-1,-1) node[below]{$\textcolor[rgb]{1.00,0.00,0.00}{J_{4}^{(1)}}$};
\end{tikzpicture}}
\centerline{\noindent {\small \textbf{Figure 4.} (Color online) The  contour  $\Sigma^{(1)}$.}}
\subsection{Second deformation}
Through introducing a scale RH problem
\begin{align}\label{43}
\left\{
\begin{array}{lr}
\delta(k)\ \mbox{is analytic in} \ \mathbb{C}\setminus (k_{1}, k_{0})\cup(k_{2}, \infty),\\
\delta_{+}(k)=\delta_{-}(k)[1+\gamma(k)\gamma^{\ast}(k)],\qquad k\in (k_{1}, k_{0})\cup(k_{2}, \infty),\\
\delta(k)\rightarrow 1,\qquad k\rightarrow \infty,
  \end{array}
\right.
\end{align}
we can delete the jump across the cut $(-\infty, k_{0})\cup(k_{2}, \infty)$.
The above RH problem \eqref{43} can be solved by Plemelj formula, given by
\begin{align}\label{45}
\delta(k)=\exp\big\{\frac{1}{2\pi i}\int_{(k_{1}, k_{0})\cup(k_{2}, \infty)}\frac{\ln[1+\gamma(y)\gamma^{\ast}(y)]}{y-k}dy\big\}.
\end{align}
To finish the second deformation, we choose the transformation
\begin{align}\label{46}
m^{(2)}=m^{(1)}\delta^{-\sigma_{3}}
\end{align}
to get a new matrix-value function $m^{(2)}$, which meets the following RH problem with the contour $\Sigma^{(2)}$ displayed in Fig. 5
\begin{align}\label{46.1}
\left\{
\begin{array}{lr}
m^{(2)}(x, t, k)\ \mbox{is analytic in} \ \mathbb{C}\setminus\Sigma^{(2)},\\
m^{(2)}_{+}(x, t, k)=m^{(2)}_{-}(x, t, k)J^{(2)}(x, t, k), \qquad k\in\Sigma^{(2)},\\
m^{(2)}(x, t, k)\rightarrow I,\qquad k\rightarrow \infty,
  \end{array}
\right.
\end{align}
and with the help of $J^{(2)}=\delta_{-}^{\sigma_{3}}J^{(1)}\delta_{+}^{-\sigma_{3}}$, we have
\begin{gather}
J_{1}^{(2)}=\left(\begin{array}{cc}
   d^{-\frac{1}{2}}  &  \delta^{2}\frac{d^{\frac{1}{2}}\gamma^{\ast}e^{-2if t}}{1+\gamma\gamma^{\ast}}\\
    0&  d^{\frac{1}{2}}\\
\end{array}\right),\ J_{2}^{(2)}=\left(\begin{array}{cc}
   d^{-\frac{1}{2}}  &  0\\
    \delta^{-2}\frac{d^{\frac{1}{2}}\gamma e^{2if t}}{1+\gamma\gamma^{\ast}}&  d^{\frac{1}{2}}\\
\end{array}\right), \ J_{\eta}^{(2)}=\left(\begin{array}{cc}
   0  &   \delta^{2}\frac{iq_{+}}{q_{0}}\\
  \delta^{-2}\frac{iq^{\ast}_{+}}{q_{0}} &  0\\
\end{array}\right),
\notag\\
J_{3}^{(2)}=\left(\begin{array}{cc}
   d^{-\frac{1}{2}}  &   0\\
  \delta^{-2}d^{-\frac{1}{2}}\gamma e^{2if t}  &  d^{\frac{1}{2}}\\
\end{array}\right),\ J_{4}^{(2)}=\left(\begin{array}{cc}
   d^{-\frac{1}{2}}  &   \delta^{2}d^{-\frac{1}{2}}\gamma^{\ast} e^{-2if t}\\
    0 &  d^{\frac{1}{2}}\\
\end{array}\right). \label{47}
\end{gather}

\centerline{\begin{tikzpicture}
\draw[-][thick](0,0)--(0,1);
\draw[-][thick](0,1)--(0,1.4);
\draw[<-][thick](0,0)--(0,-1);
\draw[-][thick](0,-1)--(0,-1.4);
\draw [-,dashed,thick] (-4,0) to [out=90,in=-60] (-5,2);
\draw [-,dashed,thick] (-4,0) to [out=-90,in=60] (-5,-2);
\draw [-,dashed,thick] (3,0) to [out=90,in=-120] (4,2);
\draw [-,dashed,thick] (3,0) to [out=-90,in=120] (4,-2);
\draw [-,dashed,thick] (0,1.4) to [out=180,in=-60] (-0.5,3);
\draw [-,dashed,thick] (0,-1.4) to [out=180,in=60] (-0.5,-3);
\draw[->][thick](-6,1)--(-5,1);
\draw [-,thick] (-5,1) to [out=-5,in=145] (-4,0);
\draw[->][thick](-6,-1)--(-5,-1);
\draw [-,thick] (-5,-1) to [out=5,in=-145] (-4,0);
\draw [->,thick] (-4,0) to [out=35,in=180] (-3,1);
\draw [-,thick] (-3,1) to [out=0,in=145] (-2,0);
\draw [->,thick] (-4,0) to [out=-35,in=-180] (-3,-1);
\draw [-,thick] (-3,-1) to [out=0,in=-145] (-2,0);
\draw [->,thick] (3,0) to [out=35,in=180] (4,1);
\draw[-][thick](4,1)--(6,1);
\draw [->,thick] (3,0) to [out=-35,in=-180] (4,-1);
\draw[-][thick](4,-1)--(6,-1);
\draw[fill] (3,0) node[below]{$k_{2}$};
\draw[fill] (-2,0) node[below]{$k_{0}$};
\draw[fill] (-4,0) node[below]{$k_{1}$};
\draw[fill] (-6,1) node[above]{$\textcolor[rgb]{1.00,0.00,0.00}{J_{3}^{(2)}}$};
\draw[fill] (-6,-1) node[below]{$\textcolor[rgb]{1.00,0.00,0.00}{J_{4}^{(2)}}$};
\draw[fill] (-3,1) node[above]{$\textcolor[rgb]{1.00,0.00,0.00}{J_{1}^{(2)}}$};
\draw[fill] (-3,-1) node[below]{$\textcolor[rgb]{1.00,0.00,0.00}{J_{2}^{(2)}}$};
\draw[fill] (5,1) node[above]{$\textcolor[rgb]{1.00,0.00,0.00}{J_{1}^{(2)}}$};
\draw[fill] (5,-1) node[below]{$\textcolor[rgb]{1.00,0.00,0.00}{J_{2}^{(2)}}$};
\draw[fill] (0.3,0.0) node{$\textcolor[rgb]{1.00,0.00,0.00}{J_{\eta}^{(2)}}$};
\draw [->,thick] (-2,0) to [out=35,in=180] (0.5,2);
\draw [-,thick] (0.5,2) to [out=0,in=145] (3,0);
\draw [->,thick] (-2,0) to [out=-35,in=-180] (0.5,-2);
\draw [-,thick] (0.5,-2) to [out=0,in=-145] (3,0);
\draw[fill] (-1,1) node[above]{$\textcolor[rgb]{1.00,0.00,0.00}{J_{3}^{(2)}}$};
\draw[fill] (-1,-1) node[below]{$\textcolor[rgb]{1.00,0.00,0.00}{J_{4}^{(2)}}$};
\draw[fill] (-3,0) node{$\widehat{\Omega}_{3}$};
\draw[fill] (-1,0) node{$\widehat{\Omega}_{4}$};
\draw[fill] (5,0) node{$\widehat{\Omega}_{3}$};
\draw[fill] (-5,0) node{$\widehat{\Omega}_{4}$};
\draw[fill] (-2,1.5) node{$\widehat{\Omega}_{1}$};
\draw[fill] (-2,-1.5) node{$\widehat{\Omega}_{2}$};
\end{tikzpicture}}
\centerline{\noindent {\small \textbf{Figure 5.} (Color online) The  contour  $\Sigma^{(2)}$.}}
\subsection{Third deformation}
For the third deformation, we select the following transformation
\begin{align}\label{48}
m^{(3)}=m^{(2)}\widehat{B}(k),
\end{align}
with
\begin{align}\label{49}
\widehat{B}(k)=\left\{
\begin{array}{lr}
d^{\frac{\sigma_{3}}{2}}\ \mbox{on} \ k\in\widehat{\Omega}_{1},\\
d^{-\frac{\sigma_{3}}{2}}\ \mbox{on} \ k\in\widehat{\Omega}_{2},\\
I \ \mbox{on} \ k\in \widehat{\Omega}_{3}\cup\widehat{\Omega}_{4},
  \end{array}
\right.
\end{align}
The goal of this transformation is to wipe out the term $\Delta(k)$.
Then the following RH problem about $m^{(3)}$ is obtained
\begin{align}\label{49.1}
\left\{
\begin{array}{lr}
m^{(3)}(x, t, k)\ \mbox{is analytic in} \ \mathbb{C}\setminus\Sigma^{(3)},\\
m^{(3)}_{+}(x, t, k)=m^{(3)}_{-}(x, t, k)J^{(3)}(x, t, k), \qquad k\in\Sigma^{(3)},\\
m^{(3)}(x, t, k)\rightarrow I,\qquad k\rightarrow \infty,
  \end{array}
\right.
\end{align}
of which the contour $\Sigma^{(3)}=\Sigma^{(2)}$ is shown in Figure 5, and $J^{(3)}$ is
\begin{gather}
J_{1}^{(3)}=\left(\begin{array}{cc}
   1  &  \delta^{2}\frac{\gamma^{\ast}e^{-2if t}}{1+\gamma\gamma^{\ast}}\\
    0&  1\\
\end{array}\right),\ J_{2}^{(3)}=\left(\begin{array}{cc}
   1  &  0\\
    \delta^{-2}\frac{\gamma e^{2if t}}{1+\gamma\gamma^{\ast}}&  1\\
\end{array}\right),
\notag\\
J_{3}^{(3)}=\left(\begin{array}{cc}
   1  &   0\\
  \delta^{-2}\gamma e^{2if t}  &  1\\
\end{array}\right),\ J_{4}^{(3)}=\left(\begin{array}{cc}
   1  &   \delta^{2}\gamma^{\ast} e^{-2if t}\\
    0 &  1\\
\end{array}\right),\ J_{\eta}^{(3)}=J_{\eta}^{(2)}. \label{50}
\end{gather}

\subsection{Eliminating of the exponential growth}
According to the sign structure of $Re(if)$ shown in Figure 2,  we see that the jump matrices $J_{3}^{(3)}$ and $J_{4}^{(3)}$ in Eq. \eqref{50} should be grow exponentially in the segment $[k_{0},\chi]$ and $[k_{0},\chi^{\ast}]$, respectively. Thus, the matrices $J_{3}^{(3)}$ and $J_{4}^{(3)}$  must be decomposed(see Figure 6)
\begin{align}\label{67}
J_{3}^{(3)}=J_{5}^{(3)}J_{7}^{(3)}J_{5}^{(3)},\qquad J_{4}^{(3)}=J_{6}^{(3)}J_{8}^{(3)}J_{6}^{(3)},
\end{align}
where
\begin{gather}\label{68}
J_{5}^{(3)}=\left(\begin{array}{cc}
   1  &  \delta^{2}\gamma^{-1}e^{-2if t}\\
   0 &  1\\
\end{array}\right),\quad J_{6}^{(3)}=\left(\begin{array}{cc}
   1  &  0\\
   \delta^{-2}(\gamma^{\ast})^{-1}e^{2if t} &  1\\
\end{array}\right),\notag\\
J_{7}^{(3)}=\left(\begin{array}{cc}
   0  &  -\delta^{2}\gamma^{-1}e^{-2if t}\\
  \delta^{-2}\gamma e^{2if t} &  0\\
\end{array}\right),\quad J_{8}^{(3)}=\left(\begin{array}{cc}
   0  &  \delta^{2}\gamma^{\ast}e^{-2if t}\\
   -\delta^{-2}(\gamma^{\ast})^{-1}e^{2if t} &  0\\
\end{array}\right).
\end{gather}

\centerline{\begin{tikzpicture}[scale=1.0]
\draw[-][thick](0,0)--(0,1);
\draw[-][thick](0,1)--(0,1.4);
\draw[<-][thick](0,0)--(0,-1);
\draw[-][thick](0,-1)--(0,-1.4);
\draw [-,dashed,thick] (-4,0) to [out=90,in=-60] (-5,2);
\draw [-,dashed,thick] (-4,0) to [out=-90,in=60] (-5,-2);
\draw [-,dashed,thick] (3,0) to [out=90,in=-120] (4,2);
\draw [-,dashed,thick] (3,0) to [out=-90,in=120] (4,-2);
\draw [-,dashed,thick] (0,1.4) to [out=180,in=-60] (-0.5,3);
\draw [-,dashed,thick] (0,-1.4) to [out=180,in=60] (-0.5,-3);
\draw[->][thick](-6,1)--(-5,1);
\draw [-,thick] (-5,1) to [out=-5,in=145] (-4,0);
\draw[->][thick](-6,-1)--(-5,-1);
\draw [-,thick] (-5,-1) to [out=5,in=-145] (-4,0);
\draw [->,thick] (-4,0) to [out=35,in=180] (-3,1);
\draw [-,thick] (-3,1) to [out=0,in=145] (-2,0);
\draw [->,thick] (-4,0) to [out=-35,in=-180] (-3,-1);
\draw [-,thick] (-3,-1) to [out=0,in=-145] (-2,0);
\draw [->,thick] (3,0) to [out=35,in=180] (4,1);
\draw[-][thick](4,1)--(6,1);
\draw [->,thick] (3,0) to [out=-35,in=-180] (4,-1);
\draw[-][thick](4,-1)--(6,-1);
\draw[fill] (3,0) node[below]{$k_{2}$};
\draw[fill] (-2,0) node[below]{$k_{0}$};
\draw[fill] (-4,0) node[below]{$k_{1}$};
\draw[fill] (-6,1) node[above]{$\textcolor[rgb]{1.00,0.00,0.00}{J_{3}^{(3)}}$};
\draw[fill] (-6,-1) node[below]{$\textcolor[rgb]{1.00,0.00,0.00}{J_{4}^{(3)}}$};
\draw[fill] (-3,1) node[above]{$\textcolor[rgb]{1.00,0.00,0.00}{J_{1}^{(3)}}$};
\draw[fill] (-3,-1) node[below]{$\textcolor[rgb]{1.00,0.00,0.00}{J_{2}^{(3)}}$};
\draw[fill] (5,1) node[above]{$\textcolor[rgb]{1.00,0.00,0.00}{J_{1}^{(3)}}$};
\draw[fill] (5,-1) node[below]{$\textcolor[rgb]{1.00,0.00,0.00}{J_{2}^{(3)}}$};
\draw[fill] (0.3,0.0) node{$\textcolor[rgb]{1.00,0.00,0.00}{J_{\eta}^{(3)}}$};
\draw [->,thick] (-2,0) to [out=35,in=180] (0.5,2);
\draw [-,thick] (0.5,2) to [out=0,in=145] (3,0);
\draw [->,thick] (-2,0) to [out=-35,in=-180] (0.5,-2);
\draw [-,thick] (0.5,-2) to [out=0,in=-145] (3,0);
\draw[fill] (1.5,1) node[above]{$\textcolor[rgb]{1.00,0.00,0.00}{J_{3}^{(3)}}$};
\draw[fill] (1.5,-1) node[below]{$\textcolor[rgb]{1.00,0.00,0.00}{J_{4}^{(3)}}$};
\draw [-,thick] (-0.35,1.72) to [out=-90,in=60] (-0.6,0.4);
\draw [->,thick] (-2,0) to [out=10,in=-120] (-0.6,0.4);
\draw [-,thick] (-0.35,-1.72) to [out=90,in=-60] (-0.6,-0.4);
\draw [->,thick] (-2,0) to [out=-10,in=120] (-0.6,-0.4);
\draw [->,thick] (-2,0) to [out=90,in=-120] (-1.75,1.35);
\draw [-,thick] (-1.75,1.35) to [out=60,in=180] (-0.35,1.72);
\draw [->,thick] (-2,0) to [out=-90,in=120] (-1.75,-1.35);
\draw [-,thick] (-1.75,-1.35) to [out=-60,in=-180] (-0.35,-1.72);
\draw[fill] (-0.35,1.72) node[right]{$\chi$};
\draw[fill] (-0.35,-1.72) node[right]{$\chi^{\ast}$};
\draw[fill] (-1.75,-1.35) node[below]{$\textcolor[rgb]{1.00,0.00,0.00}{J_{6}^{(3)}}$};
\draw[fill] (-1,-0.7) node[below]{$\textcolor[rgb]{1.00,0.00,0.00}{J_{8}^{(3)}}$};
\draw[fill] (-0.4,-0.9) node[above]{$\textcolor[rgb]{1.00,0.00,0.00}{J_{6}^{(3)}}$};
\draw[fill] (-2,1.4) node[below]{$\textcolor[rgb]{1.00,0.00,0.00}{J_{5}^{(3)}}$};
\draw[fill] (-1.1,1) node[below]{$\textcolor[rgb]{1.00,0.00,0.00}{J_{7}^{(3)}}$};
\draw[fill] (-0.4,0.9) node[below]{$\textcolor[rgb]{1.00,0.00,0.00}{J_{5}^{(3)}}$};
\end{tikzpicture}}
\centerline{\noindent {\small \textbf{Figure 6.} (Color online) The contour  $\widehat{\Sigma}^{(3)}$.}}

Next, we set the  transformation $m^{(4)}=m^{(3)}e^{iG(k)t\sigma_{3}}$ via using a time-dependent $G$ function which is analytic off the cuts $\eta\cup\varpi$. Of which $\varpi=\varpi_{+}\cup\varpi_{-}$ with $\varpi_{+}=[k_{0},\chi]$ and
$\varpi_{-}=[k_{0},\chi^{\ast}]$. Then the new jump matrices $J^{(4)}$ is calculated as follows
\begin{gather}
J_{1}^{(4)}=\left(\begin{array}{cc}
   1  &  \delta^{2}\frac{\gamma^{\ast}e^{-2i(f+G)t}}{1+\gamma\gamma^{\ast}}\\
    0&  1\\
\end{array}\right),\ J_{2}^{(4)}=\left(\begin{array}{cc}
   1  &  0\\
    \delta^{-2}\frac{\gamma e^{2i(f+G)t}}{1+\gamma\gamma^{\ast}}&  1\\
\end{array}\right),
\notag\\
J_{3}^{(4)}=\left(\begin{array}{cc}
   1  &   0\\
  \delta^{-2}\gamma e^{2i(f+G)t}  &  1\\
\end{array}\right),\ J_{4}^{(4)}=\left(\begin{array}{cc}
   1  &   \delta^{2}\gamma^{\ast} e^{-2i(f+G)t}\\
    0 &  1\\
\end{array}\right),\notag\\
\ J_{\eta}^{(4)}=\left(\begin{array}{cc}
   0  &   \delta^{2}\frac{iq_{+}}{q_{0}}e^{-i(G_{+}+G_{-})t}\\
  \delta^{-2}\frac{iq^{\ast}_{+}}{q_{0}}e^{i(G_{+}+G_{-})t} &  0\\
\end{array}\right),\notag\\
J_{5}^{(4)}=\left(\begin{array}{cc}
   1  &  \delta^{2}\gamma^{-1}e^{-2i(f+G)t}\\
   0 &  1\\
\end{array}\right),\quad J_{6}^{(4)}=\left(\begin{array}{cc}
   1  &  0\\
   \delta^{-2}(\gamma^{\ast})^{-1}e^{2i(f+G)t} &  1\\
\end{array}\right),\notag\\
J_{7}^{(4)}=\left(\begin{array}{cc}
   0  &  -\delta^{2}\gamma^{-1}e^{-i(2f+G_{+}+G_{-})t}\\
  \delta^{-2}\gamma e^{i(2f+G_{+}+G_{-})t} &  0\\
\end{array}\right),\notag\\ J_{8}^{(4)}=\left(\begin{array}{cc}
   0  &  \delta^{2}\gamma^{\ast}e^{-i(2f+G_{+}+G_{-})t}\\
   -\delta^{-2}(\gamma^{\ast})^{-1}e^{i(2f+G_{+}+G_{-})t} &  0\\
\end{array}\right). \label{70}
\end{gather}
Furthermore, we introduce a new function $\omega$:
\begin{align}\label{71}
\omega(k)=f(k)+G(k),
\end{align}
whose  properties need to be investigated up to find the parameters $k_{0}$ and  $\chi$. First of all, a function $z$ should be defined as
\begin{align}\label{72}
z(k)=\sqrt{(k^{2}+q_{0}^{2})(k-\chi)(k-\chi^{\ast})},
\end{align}
which has branch cuts $\eta\cup\varpi$ and satisfies $z(k)=-z_{+}(k)=z_{-}(k)$. We implement this algebraic curve as two Riemann surfaces, and  the basis $\{L_{1}, L_{2}\}$ cycles of this Riemann surface can be described as: the $L_{1}$-cycle is a simple counterclockwise closed ring around the bifurcation incision $\eta$, which lies on the lower
sheet. The $L_{2}$-cycle starts  from the point $\chi$ on the upper sheet, then accesses $-iq_{0}$ and
gets back to the starting point on the lower sheet.

Next, we let $\omega(k)$ meets
\begin{align}\label{73}
\omega(k)=\frac{1}{2}\left(\int_{iq_{0}}^{k}
+\int_{-iq_{0}}^{k}\right)d\omega(y),
\end{align}
which is a Abelian integral and $d\omega$ is given by
\begin{align}\label{74}
d\omega(k)=-32\frac{(k-k_{0})(k-k_{1})(k-k_{2})(k-\chi)(k-\chi^{\ast})}{z(k)}dk.
\end{align}
Moreover, the sign signatures of Im$\omega(k)$ must be same as that ones of Im$f(k)$ for large $k$, given by
\begin{align}\label{78}
\mbox{Im} \omega=\mbox{Im} f+\mathcal{O}(\frac{1}{k}),\qquad k\rightarrow\infty.
\end{align}
Therefore, we have
\begin{align}\label{85}
\omega(k)=-8k^{4}+2k^{2}+\xi k+\omega_{0}+\mathcal{O}(\frac{1}{k}), \quad k\rightarrow\infty.
\end{align}
When $k\rightarrow\infty$, the large $k$ expression of $z(k)$ become
\begin{align}\label{79}
z(k)=k^{2}\left[1-\frac{\chi+\chi^{\ast}}{2k}
+\frac{4q_{0}^{2}-(\chi-\chi^{\ast})^{2}}{8k^{2}}+\mathcal{O}(\frac{1}{k^{3}})\right],\quad k\rightarrow\infty.
\end{align}
Taking $\chi=\chi_{1}+\chi_{2}i$, from \eqref{74}, we easily obtain
\begin{align}\label{80}
\frac{d\omega}{dk}=&-32k^{3}+32(\chi_{1}+k_{0}+k_{1}+k_{2})k^{2}-[32(k_{0}+k_{1}+k_{2})\chi_{1}
+32(k_{1}k_{0}+k_{2}k_{0}+k_{2}k_{1})\notag\\
&-16q_{0}^{2}+16\chi_{2}^{2}]k
+(32(k_{1}k_{0}+k_{2}k_{0}+k_{2}k_{1})-16q_{0}^{2}+16\chi_{2}^{2})\chi_{1}+32k_{0}k_{1}k_{2}\notag\\
&
+16(\chi_{2}^{2}-q_{0}^{2})(k_{0}+k_{1}+k_{2})+\mathcal{O}(\frac{1}{k}), \quad k\rightarrow\infty.
\end{align}
Meanwhile, since $f=\lambda\left[\xi-8k^{3}+2k+4kq_{0}^{2}\right]$, one has
\begin{align}\label{82}
\frac{df(k)}{dk}=-32k^{3}+4k+\xi+\mathcal{O}(\frac{1}{k}), \quad k\rightarrow\infty.
\end{align}
Since \eqref{78} is allowed, and from Eqs. \eqref{80}, \eqref{82}, we can derive
\begin{align}\label{84}
\chi_{1}=&-k_{0}-k_{1}-k_{2},\notag\\ \chi_{2}=&\frac{1}{2}\sqrt{8(k_{0}^{2}+k_{1}^{2}+k_{2}^{2})
+8(k_{1}k_{0}+k_{2}k_{0}+k_{2}k_{1})+4q_{0}^{2}-1},
\end{align}
where parameter $k_{0}$ is still need to be derived later.
Observing that
\begin{align}\label{85}
-16\left(\int_{iq_{0}}^{k}
+\int_{-iq_{0}}^{k}\right)(y^{3}-\frac{y}{8}-\frac{\xi}{32})dy
=-8k^{4}+2k^{2}+\xi k+8q_{0}^{4}+2q_{0}^{2},
\end{align}
then the expression of $\omega(k)$ in Eq.\eqref{73} can been reconstructed as
\begin{align}\label{86}
\omega(k)=&-16\left(\int_{iq_{0}}^{k}
+\int_{-iq_{0}}^{k}\right)[\frac{(y-k_{0})(y-k_{1})(y-k_{2})(y-\chi)(y-\chi^{\ast})}{z(y)}\notag\\
&-(y^{3}-\frac{y}{8}-\frac{\xi}{32})]dy-8k^{4}+2k^{2}+\xi k+8q_{0}^{4}+2q_{0}^{2}.
\end{align}
As $k\rightarrow\infty$ in \eqref{86}, we have
\begin{align}\label{87}
\omega_{0}=&-16\left(\int_{iq_{0}}^{\infty}
+\int_{-iq_{0}}^{\infty}\right)[\frac{(y-k_{0})(y-k_{1})(y-k_{2})(y-\chi)(y-\chi^{\ast})}{z(y)}\notag\\
&-(y^{3}-\frac{y}{8}-\frac{\xi}{32})]dy+8q_{0}^{4}+2q_{0}^{2}.
\end{align}

Next, we will devote to reveal the parameter $k_{0}$ on the real line by presenting the asymptotic expansions of $\omega(k)$ near point $\chi$.
Similar to reference \cite{Liunan7}, it is not hard to obtain that
\begin{align}\label{94}
&\int_{-iq_{0}}^{iq_{0}}[\frac{(y-k_{0})(y-k_{1})(y-k_{2})(y-\chi)(y-\chi^{\ast})}{z(y)}]dy\notag\\
&=\int_{-iq_{0}}^{iq_{0}}\sqrt{\frac{(y-\chi_{1})^{2}+\chi_{2}^{2}}{y^{2}+q_{0}^{2}}}(y-k_{0})(y-k_{1})(y-k_{2})dy=0,
\end{align}
which uniquely gives us the point $k_{0}$ is uniquely expressed.

Now, the function $\omega(k)$ satisfies the following jump condition:
\begin{align}\label{95}
&\omega_{+}(k)+\omega_{-}(k)=0,\qquad k\in \eta,\notag\\
&\omega_{+}(k)+\omega_{-}(k)=\Omega,\qquad k\in \varpi,
\end{align}
where  $\Omega$ is real constant given by
\begin{align}\label{96}
\Omega=-32\left(\int_{iq_{0}}^{\chi}
+\int_{-iq_{0}}^{\chi^{\ast}}\right)\frac{(k-k_{0})(k-k_{1})(k-k_{2})(k-\chi)(k-\chi^{\ast})}{z(k)}dk.
\end{align}
Besides, since function $\omega(k)$ is defined in Eq. \eqref{71}, one easily obtains
\begin{align}\label{99}
G(\infty)=\omega_{0}-3q_{0}^{4}-q_{0}^{2}, \quad k\rightarrow\infty,
\end{align}
and we also have $m^{(4)}\rightarrow e^{iG(\infty)t\sigma_{3}}$ as $k\rightarrow\infty$.
Finally, we can derive the RH problem for $m^{(4)}$,
whose jump matrices $J^{(4)}$ are
\begin{gather}
J_{1}^{(4)}=\left(\begin{array}{cc}
   1  &  \delta^{2}\frac{\gamma^{\ast}e^{-2i\omega t}}{1+\gamma\gamma^{\ast}}\\
    0&  1\\
\end{array}\right),\ J_{2}^{(4)}=\left(\begin{array}{cc}
   1  &  0\\
    \delta^{-2}\frac{\gamma e^{2i\omega t}}{1+\gamma\gamma^{\ast}}&  1\\
\end{array}\right),
\notag\\
J_{3}^{(4)}=\left(\begin{array}{cc}
   1  &   0\\
  \delta^{-2}\gamma e^{2i\omega t}  &  1\\
\end{array}\right),\ J_{4}^{(4)}=\left(\begin{array}{cc}
   1  &   \delta^{2}\gamma^{\ast} e^{-2i\omega t}\\
    0 &  1\\
\end{array}\right),\notag\\
\ J_{\eta}^{(4)}=\left(\begin{array}{cc}
   0  &   \delta^{2}\frac{iq_{+}}{q_{0}}\\
  \delta^{-2}\frac{iq^{\ast}_{+}}{q_{0}} &  0\\
\end{array}\right),\notag\\
J_{5}^{(4)}=\left(\begin{array}{cc}
   1  &  \delta^{2}\gamma^{-1}e^{-2i\omega t}\\
   0 &  1\\
\end{array}\right),\quad J_{6}^{(4)}=\left(\begin{array}{cc}
   1  &  0\\
   \delta^{-2}(\gamma^{\ast})^{-1}e^{2i\omega t} &  1\\
\end{array}\right),\notag\\
J_{7}^{(4)}=\left(\begin{array}{cc}
   0  &  -\delta^{2}\gamma^{-1}e^{-i\Omega t}\\
  \delta^{-2}\gamma e^{i\Omega t} &  0\\
\end{array}\right),\notag\\ J_{8}^{(4)}=\left(\begin{array}{cc}
   0  &  \delta^{2}\gamma^{\ast}e^{-i\Omega t}\\
   -\delta^{-2}(\gamma^{\ast})^{-1}e^{i\Omega t} &  0\\
\end{array}\right). \label{100}
\end{gather}
The sign signature of Im$(\omega)(k)$ ensures that the jump matrices $J_{i}^{(4)}(i = 1, 2, 3, 4, 5, 6)$ are all exponentially decaying in the associated branch cuts.
\subsection{Further deformation}
In order to delete the variable $k$ from the jump matrices  $J_{\eta}^{(4)}, J_{7}^{(4)}, J_{8}^{(4)}$, we need to introduce the $g$-function mechanism again. In the same way, we select following transformation
\begin{align}\label{101}
m^{(5)}=m^{(4)}g(k)^{\sigma_{3}},
\end{align}
where the function $g(k)$, which is analytic in $\mathbb{C}\setminus (\eta\cup\varpi)$, satisfies
\begin{align}\label{102}
g_{+}(k)g_{-}(k)=\left\{
\begin{array}{lr}
\delta^{2}\qquad\qquad\ \mbox{on} \qquad k\in\eta,\\
e^{i\vartheta}\frac{\delta^{2}}{\gamma}\qquad\quad \mbox{on} \qquad k\in\varpi_{+},\\
e^{i\vartheta}\delta^{2}\gamma^{\ast}\qquad \mbox{on} \qquad k\in\varpi_{-},
  \end{array}
\right.
\end{align}
of which $\vartheta$ is a real constant and given by
\begin{align}\label{103a}
\vartheta=i\frac{\int_{\eta}\frac{2\ln\delta}{z}ds+\int_{[k_{0},\chi]}\frac{\ln\frac{\delta^{2}}{\gamma}}{z}ds
+\int_{[k_{0},\chi^{\ast}]}\frac{\ln\delta^{2}\gamma^{\ast}}{z}ds}{\int_{[k_{0},\chi]\cup[k_{0},\chi^{\ast}]}\frac{1}{z}ds}.
\end{align}
Applying the Plemelj's formula, the $g(k)$ function can be solved by the following integral representation
\begin{align}\label{103b}
g(k)=\exp\{-\frac{z}{2\pi i}(\int_{\eta}\frac{2\ln\delta}{z(s-k)}ds+\int_{[k_{0},\chi]}\frac{i\vartheta+\ln\frac{\delta^{2}}{\gamma}}{z(s-k)}ds
+\int_{[k_{0},\chi^{\ast}]}\frac{i\vartheta+\ln\delta^{2}\gamma^{\ast}}{z(s-k)}ds)\},
\end{align}
which implies that $g(k)$ has the following behavior for the large $k$:
\begin{align}\label{103c}
g(k)=e^{ig(\infty)}+\mathcal{O}(\frac{1}{k}), \quad k\rightarrow\infty,
\end{align}
where the $g(\infty)$ is a real constant, given by
\begin{align}\label{103}
g(\infty)=-\frac{1}{2\pi}(\int_{\eta}\frac{2\ln\delta}{z}sds+\int_{[k_{0},\chi]}\frac{i\vartheta+\ln\frac{\delta^{2}}{\gamma}}{z}sds
+\int_{[k_{0},\chi^{\ast}]}\frac{i\vartheta+\ln\delta^{2}\gamma^{\ast}}{z}sds).
\end{align}

Finally, we get the following RH problem for $m^{(5)}$
\begin{align}\label{104a}
\left\{
\begin{array}{lr}
m^{(5)}(x, t, k)\ \mbox{is analytic in} \ \mathbb{C}\setminus\Sigma^{(5)},\\
m^{(5)}_{+}(x, t, k)=m^{(5)}_{-}(x, t, k)J^{(5)}(x, t, k), \qquad k\in\Sigma^{(5)},\\
m^{(5)}(x, t, k)\rightarrow e^{i(g(\infty)+G(\infty)t)\sigma_{3}},\qquad k\rightarrow \infty,
  \end{array}
\right.
\end{align}
where the contour $\Sigma^{(5)}=\widehat{\Sigma}^{(3)}$ and the jump matrices $J^{(5)}$ become
\begin{gather}
J_{1}^{(5)}=\left(\begin{array}{cc}
   1  &  \delta^{2}\frac{\gamma^{\ast}e^{-2i\omega t}g^{-2}}{1+\gamma\gamma^{\ast}}\\
    0&  1\\
\end{array}\right),\ J_{2}^{(5)}=\left(\begin{array}{cc}
   1  &  0\\
    \delta^{-2}\frac{\gamma e^{2i\omega t}g^{2}}{1+\gamma\gamma^{\ast}}&  1\\
\end{array}\right),
\notag\\
J_{3}^{(5)}=\left(\begin{array}{cc}
   1  &   0\\
  \delta^{-2}\gamma e^{2i\omega t}g^{2}  &  1\\
\end{array}\right),\ J_{4}^{(5)}=\left(\begin{array}{cc}
   1  &   \delta^{2}\gamma^{\ast} e^{-2i\omega t}g^{-2}\\
    0 &  1\\
\end{array}\right),\notag\\
\ J_{\eta}^{(5)}=\left(\begin{array}{cc}
   0  &   \frac{iq_{+}}{q_{0}}\\
  \frac{iq^{\ast}_{+}}{q_{0}} &  0\\
\end{array}\right),\notag\\
J_{5}^{(5)}=\left(\begin{array}{cc}
   1  &  \delta^{2}\gamma^{-1}e^{-2i\omega t}g^{-2}\\
   0 &  1\\
\end{array}\right),\quad J_{6}^{(5)}=\left(\begin{array}{cc}
   1  &  0\\
   \delta^{-2}(\gamma^{\ast})^{-1}e^{2i\omega t}g^{2} &  1\\
\end{array}\right),\notag\\
J_{7}^{(5)}=\left(\begin{array}{cc}
   0  &  -e^{-i(\Omega t+\vartheta)}\\
   e^{i(\Omega t+\vartheta)} &  0\\
\end{array}\right),\notag\\ J_{8}^{(5)}=\left(\begin{array}{cc}
   0  &  e^{-i(\Omega t+\vartheta)}\\
  -e^{i(\Omega t+\vartheta)} &  0\\
\end{array}\right). \label{104}
\end{gather}
\subsection{Model problem and the results}
Since the jump matrices $J_{i}^{(5)}(i = 1, 2, 3, 4, 5, 6)$ are
all exponentially decaying to the identity away from the points $k_{0}, k_{1}, k_{2}$, $\chi$ and $\chi^{\ast}$ as $t\rightarrow\infty$, we can obtain a  model problem to determine the leading term of the solution, given by
\begin{align}\label{105}
\left\{
\begin{array}{lr}
m^{mod}(x, t, k)\ \mbox{is analytic in} \ \mathbb{C}\setminus (\eta\cup\varpi_{+}\cup(-\varpi_{-})),\\
m^{mod}_{+}(x, t, k)=m^{mod}_{-}(x, t, k)J^{mod}(x, t, k), \qquad k\in\eta\cup\varpi_{+}\cup(-\varpi_{-}),\\
m^{mod}(x, t, k)\rightarrow e^{i(g(\infty)+G(\infty)t)\sigma_{3}},\qquad k\rightarrow \infty,
  \end{array}
\right.
\end{align}
where
\begin{gather}
J_{\eta}^{mod}=J_{\eta}^{(5)}=\left(\begin{array}{cc}
   0  &   \frac{iq_{+}}{q_{0}}\\
  \frac{iq^{\ast}_{+}}{q_{0}} &  0\\
\end{array}\right),\notag\\
J_{\varpi_{+}\cup(-\varpi_{-})}^{mod}=\left(\begin{array}{cc}
   0  &  -e^{-i(\Omega t+\vartheta)}\\
   e^{i(\Omega t+\vartheta)} &  0\\
\end{array}\right), \label{106}
\end{gather}
and $-\varpi_{-}$ means the negative direction of cut $\varpi_{-}$.

For large $k$, introducing the factorization $m^{(5)}=m^{err}m^{mod}$ and taking the Laurent series for matrices $m^{err}, m^{mod}$ as
\begin{align}\label{106a}
&m^{err}=I+\frac{m^{err}_{1}(x, t)}{k}+\frac{m^{err}_{2}(x, t)}{k^{2}}+\mathcal{O}(\frac{1}{k^{3}}),\qquad k\rightarrow \infty,\notag\\
&m^{mod}=e^{i(g(\infty)+G(\infty)t)\sigma_{3}}+\frac{m^{mod}_{1}(x, t)}{k}+\frac{m^{mod}_{2}(x, t)}{k^{2}}+\mathcal{O}(\frac{1}{k^{3}}),\qquad k\rightarrow \infty,
\end{align}
we can represent the solution $q(x,t)$ for the Eq.\eqref{1} via the solution of model problem
\begin{align}\label{107}
q(x,t)=2i\left(m_{1}^{mod}(x,t)e^{i(g(\infty)+G(\infty)t)}+m_{1}^{err}(x,t)\right)_{12}.
\end{align}
Similar to reference \cite{Liunan7}, we get $\mid m_{1}^{err}\mid=\mathcal{O}(t^{-\frac{1}{2}})$.

To solve the model RH problem \eqref{105}, we first
define the Abelian differential in what follows
\begin{align}\label{108}
d\vartheta=\frac{\vartheta_{0}}{z(k)}dk,\qquad \vartheta_{0}=\left(\oint_{L_{1}}\frac{1}{z(k)}dk\right)^{-1},
\end{align}
which is normalized at the case of $\oint_{L_{1}}d\vartheta=1$. At the same time, the above Abelian differential \eqref{108} admits following Riemann period $\tau$
\begin{align}\label{109}
\tau=\oint_{L_{2}}d\vartheta,
\end{align}
which is purely imaginary when $i\tau<0$\cite{Farkas}. Therefore,  the theta function can be written into
\begin{align}\label{114}
\Theta(k)=\sum_{\varrho\in\mathbb{Z}}e^{2\pi i\varrho k+\pi i\tau\varrho^{2}},
\end{align}
which yields the properties
\begin{align}\label{115}
\Theta(k+n)=\Theta(k),\quad \Theta(k+n\tau)=e^{-(2\pi in k+\pi i\tau n^{2})}\Theta(k),\quad  n\in\mathbb{Z}.
\end{align}
According to the Abelian map
\begin{align}\label{110}
V(k)=\int_{iq_{0}}^{k}d\vartheta,
\end{align}
it arrives at
\begin{align}\label{111}
&V_{+}(k)+V_{-}(k)=n-\tau, \quad n\in\mathbb{Z},\quad k\in\varpi_{+}\cup(-\varpi_{-}),\notag\\
&V_{+}(k)+V_{-}(k)=n, \quad n\in\mathbb{Z},\quad k\in\eta.
\end{align}

Finally, a $2\times 2$ matrix-valued function $M(k)=M(x, t, k)$ are constructed to solve the mod problem \eqref{105}, whose elements are
\begin{align}\label{116}
M_{11}(k)=\frac{1}{2}[r(k)+r^{-1}(k)]\frac{\Theta\left(\frac{\Omega t+\vartheta+i\ln\left(\frac{iq_{+}^{\ast}}{q_{0}}\right)}{2\pi}+V(k)+C\right)}
{\sqrt{\frac{q_{0}}{iq_{+}^{\ast}}}\Theta\left(V(k)+C\right)},\notag\\
M_{12}(k)=\frac{i}{2}[r(k)-r^{-1}(k)]\frac{\Theta\left(\frac{\Omega t+\vartheta+i\ln\left(\frac{iq_{+}^{\ast}}{q_{0}}\right)}{2\pi}-V(k)+C\right)}
{\sqrt{\frac{iq_{+}^{\ast}}{q_{0}}}\Theta\left(-V(k)+C\right)},\notag\\
M_{21}(k)=\frac{-i}{2}[r(k)-r^{-1}(k)]\frac{\Theta\left(\frac{\Omega t+\vartheta+i\ln\left(\frac{iq_{+}^{\ast}}{q_{0}}\right)}{2\pi}+V(k)-C\right)}
{\sqrt{\frac{q_{0}}{iq_{+}^{\ast}}}\Theta\left(V(k)-C\right)},\notag\\
M_{22}(k)=\frac{1}{2}[r(k)+r^{-1}(k)]\frac{\Theta\left(\frac{\Omega t+\vartheta+i\ln\left(\frac{iq_{+}^{\ast}}{q_{0}}\right)}{2\pi}-V(k)-C\right)}
{\sqrt{\frac{iq_{+}^{\ast}}{q_{0}}}\Theta\left(-V(k)-C\right)},
\end{align}
of which the function $r(k)$ is
\begin{align}\label{112}
r(k)=\left(\frac{(k-\chi)(k-iq_{0})}{(k-\chi^{\ast})(k+iq_{0})}\right)^{\frac{1}{4}},
\end{align}
which has the identical jump discontinuity across $\eta$ and $\varpi_{+}\cup(-\varpi_{-})$, as well as
$r_{+}(k)=ir_{-}(k)$, and it's large-$k$ asymptotic is
\begin{align}\label{113}
&r(k)=1-\frac{i(\chi_{2}+q_{0})}{2k}+O\left(\frac{1}{k^{2}}\right),\quad k\rightarrow\infty,\notag\\
&r(k)-r^{-1}(k)=-\frac{i(\chi_{2}+q_{0})}{k}+O\left(\frac{1}{k^{2}}\right),\quad k\rightarrow\infty.
\end{align}
Besides, we also have
\begin{align}\label{117a}
C=V(\hat{k})+\frac{1}{2}(1+\tau),\qquad \hat{k}=\frac{q_{0}\chi_{1}}{q_{0}+\chi_{2}}.
\end{align}
Then the model RH problem \eqref{105} is solved as
\begin{align}\label{117}
m^{mod}(x,t,k)=e^{i(g(\infty)+G(\infty)t)\sigma_{3}}M^{-1}(\infty, C)M(k, C),
\end{align}
further we obtain
\begin{align}\label{118}
(m_{1}^{mod})_{12}=\frac{q_{0}(\chi_{2}+q_{0})\Theta\left(\frac{\Omega t+\vartheta+i\ln\left(\frac{iq_{+}^{\ast}}{q_{0}}\right)}{2\pi}-V(\infty)+C\right)
\Theta(V(\infty)+C)}{2iq_{+}^{\ast}\Theta\left(\frac{\Omega t+\vartheta+i\ln\left(\frac{iq_{+}^{\ast}}{q_{0}}\right)}{2\pi}+V(\infty)+C\right)\Theta(-V(\infty)+C)e^{-i(g(\infty)+G(\infty)t)}},
\end{align}
which implies the long-time asymptotics of solution $q(x,t)$ for the fourth-order dispersive NLS
equation \eqref{1} is
\begin{gather}
q(x,t)=\frac{q_{0}}{q_{+}^{\ast}}(\chi_{2}+q_{0})\frac{\Theta\left(\frac{\Omega t+\vartheta+i\ln\left(\frac{iq_{+}^{\ast}}{q_{0}}\right)}{2\pi}-V(\infty)+C\right)
\Theta(V(\infty)+C)}{\Theta\left(\frac{\Omega t+\vartheta+i\ln\left(\frac{iq_{+}^{\ast}}{q_{0}}\right)}{2\pi}+V(\infty)+C\right)\Theta(-V(\infty)+C)}\notag\\e^{2i(g(\infty)+G(\infty)t)}+\mathcal{O}(t^{-\frac{1}{2}}),\label{119}
\end{gather}
where $V(\infty)=\int_{iq_{0}}^{\infty}d\vartheta$, and $\chi_{2}, \Omega, G(\infty), \vartheta, g(\infty), C$ are presented in Eqs. \eqref{84}, \eqref{96}, \eqref{99}, \eqref{103a}, \eqref{103}, \eqref{117a}, respectively.

\section*{Acknowledgements}
This work was supported by the project is supported by National Natural Science Foundation of China(No.12175069)
and Science and Technology Commission of Shanghai Municipality (No.21JC1402500 and No.18dz2271000).






\begin{thebibliography}{99}
\bibitem{Zhanghaiq}
H. Q. Zhang, B. Tian, X. H. Meng, X. L\"{u}, W. J. Liu, Conservation laws, soliton solutions and modulational instability for
the higher-order dispersive nonlinear Schr\"{o}dinger equation, Eur. Phys. J. B 72 (2009) 233-239.
\bibitem{Tian22}
T. A. Davydova, Y. A. Zaliznyak, Schr\"{o}dinger ordinary solitons and chirped solitons: fourth-order dispersive effects and cubic-quintic nonlinearity, Phys. D, 156 (2001) 260-282.
\bibitem{Tian23}
F. Azzouzi, H. Triki, K. Mezghiche, A. E. Akrmi, Solitary wave solutions for high
dispersive cubic-quintic nonlinear Schr\"{o}dinger equation , Chaos Solitons Fract.
39 (2009) 1304-1307.
\bibitem{Tian24}
S. L. Palacios, J. M. Fern\'{a}ndez-D\'{i}az, Black optical solitons for media with
parabolic nonlinearity law in the presence of fourth order dispersion, Opt. Commun. 178 (2000) 457-460.
\bibitem{Tian25}
M. Daniel, L. Kavitha, R. Amuda, Soliton spin excitations in an anisotropic
Heisenberg ferromagnet with octupole-dipole interaction, Phys. Rev. B, 59 (1999)
13774-13781.
\bibitem{Tian26}
K. Porsezian, M. Daniel, M. Lakshmanan, On the integrability aspects of the onedimensional classical continuum isotropic biquadratic Heisenberg spin chain, J. Math. Phys. 33 (1992) 1807-1816.
\bibitem{Tian28}
L. H. Wang, K. Porsezian, J. S. He, Breather and rogue wave solutions of generalized nonlinear Schr\"{o}dinger equation, Phys. Rev. E, 87 (2013) 053202.
\bibitem{Tian-Li}
Z. Q. Li, S.  F.  Tian,  J. J. Yang, Riemann-Hilbert approach and soliton solutions for the higher-order dispersive nonlinear Schr\"{o}dinger equation with nonzero boundary conditions. arXiv preprint arXiv:1911.01624, 2019.
\bibitem{guo-liu9}
R. Guo, H. Q. Hao, Breathers and multi-soliton solutions for the higher-order generalized nonlinear Schr\"{o}dinger equation,
Commun. Nonlinear Sci. Numer. Simul. 18 (2013) 2426-2435.
\bibitem{guo-liu18}
R. X. Liu, B. Tian, L. C. Liu, B. Qin, X. L\"{u}, Bilinear forms, N-soliton solutions and soliton interactions for a fourth-order
dispersive nonlinear Schr\"{o}dinger equation in condensed-matter physics and biophysics, Physica B 413 (2013) 120-125.
\bibitem{guo-liu22}
X. L. Wang, W. G. Zhang, B. G. Zhai, H. Q. Zhai, Rogue waves of the higher-order dispersive nonlinear Schr\"{o}dinger equation,
Commun. Theor. Phys. 58 (2012) 531-538.
\bibitem{guo-liu}
Y. F. Wang,  N. Liu,  B. L. Guo, Long-time asymptotic behavior for a fourth-order dispersive nonlinear Schr\"{o}dinger equation. J. Math. Anal. Appl., 506(1)(2022) 125560.
\bibitem{guo-liu11}
L. Huang, Asymptotics behavior for the integrable nonlinear Schr\"{o}dinger equation with quartic terms: Cauchy problem,
J. Nonlinear Math. Phys. 27 (2020) 592-615.
\bibitem{Peng28}
S. V. Manakov,   Nonlinear Fraunnhofer diffraction, Sov. Phys. JETP, 38 (1974)  693-696.
\bibitem{Peng29}
M. J. Ablowitz,   A. C. Newell,  The decay of the continuous spectrum for solutions of the
Korteweg-de Vries equation, J. Math. Phys., 14 (1973)  1277-1284.
\bibitem{Peng30}
V. E. Zakharov,  S. V. Manakov,  Asymptotic behavior of non-linear wave systems integrated
by the inverse scattering method, Sov. Phys. JETP, 44 (1976) 106-112.
\bibitem{Peng31}
M. J. Ablowitz,  H. Segur,  Asymptotic solutions of the Korteweg-de Vries equation, Stud.
Appl. Math., 57 (1977)  13-44.
\bibitem{Peng32}
H. Segur,  M. J. Ablowitz,  Asymptotic solutions and conservation laws for the nonlinear
Schrodinger equation I, J. Math. Phys., 17 (1973) 710-713.
\bibitem{Peng33}
A. R. Its,  Asymptotics of solutions of the nonlinear Schr\"{o}dinger equation and isomonodromic
deformations of systems of linear differential equations, Sov. Math. Dokl., 24 (1981) 452-456.
\bibitem{Peng34}
P. A. Deift, X.  Zhou,  A steepest descent method for oscillatory Riemann-Hilbert problems,
Ann. Math. 137 (1993)  295-368.
\bibitem{Peng35}
P. Deift,   S. Venakides,  X. Zhou, The collisionless shock region for the long-time behavior of
solutions of the KdV equation, Commun. Pure Appl. Math., 47 (1994) 199-206.
\bibitem{Peng36}
P. Deift,  X. Zhou,  Asymptotics for the Painleve II equation, Commun. Pure Appl. Math., 48 (1995) 277-337.
\bibitem{Peng37}
P. Deift,    S. Venakides,   X. Zhou, New results in small dispersion KdV by an extension of the
steepest descent method for Riemann-Hilbert problems, Int. Math. Res. Not., (1997)
286-299.
\bibitem{Peng38}
S. Kamvissis,  Long time behavior for the focusing nonlinear Schr\"{o}edinger equation with real
spectral singularities, Commun. Math. Phys., 180  (1996) 325-341.
\bibitem{Peng39}
K. Grunert,  G. Teschl,   Long-time asymptotics for the Korteweg-de Vries equation via nonlinear steepest descent, Math. Phys. Anal. Geom. 12 (2009)  287-324.
\bibitem{Peng40}
P. J. Cheng,  S. Venakides,  X. Zhou,  Long-time asymptotics for the pure radiation solution of
the sine-Gordon equation, Commun. Partial Differential Equations, 24 (1999) 1195-1262.
\bibitem{Peng41}
A. Boutet de Monvel,  A. Kostenko,  D. Shepelsky,  G. Teschl,  Long-time asymptotics for the
Camassa-Holm equation, SIAM J. Math. Anal. 41 (2009)  1559-1588.
\bibitem{Peng43}
B. L. Guo,  N. Liu,   Y. F. Wang, Long-time asymptotics for the Hirota equation on the half-line, Nonlinear Anal., 174(2018) 118-140.
\bibitem{Peng44}
L. Huang,  J. Xu,  E. G. Fan,  Long-time asymptotic for the Hirota equation via nonlinear steepest descent method,
Nonlinear Anal. Real World Appl. 26 (2015)  229-262.
\bibitem{Wang-ke}
D. S. Wang,   X. L. Wang, Long-time asymptotics and the bright N-soliton solutions of the Kundu-Eckhaus equation
via the Riemann-Hilbert approach, Nonlinear Anal. Real World Appl. 41 (2018) 334-361.
\bibitem{Geng}
X. G. Geng,  H. Liu,  The nonlinear steepest descent method to long-time asymptotics of the
coupled nonlinear Schr\"{o}dinger equation, J. Nonlinear Sci. 28 (2018) 739-763.
\bibitem{Peng46}
A. Boutet de Monvel, A. Its,  V. Kotlyarov,  Long-time asymptotics for the focusing NLS
equation with time-periodic boundary condition on the half-line, Commun. Math. Phys. 290,
479-522(2009).
\bibitem{Peng47}
S. F. Tian,  T. T. Zhang,  Long-time asymptotic behavior for the Gerdjikov-Ivanov type of
derivative nonlinear Schr\"{o}dinger equation with time-periodic boundary condition, Proc. Am.
Math. Soc. 146 (2018) 1713-1729.
\bibitem{Peng45}
R. Buckingham,  S. Venakides,  Long-time asymptotics of the nonlinear Schr\"{o}dinger equation
shock problem, Commun. Pure Appl. Math. 60 (2007) 1349-1414.
\bibitem{Peng48}
A. Boutet de Monvel,  V. P. Kotlyarov,  D. Shepelsky,  Focusing NLS equation: Long-time
dynamics of step-like initial data, Int. Math. Res. Not. (2011) 1613-1653.
\bibitem{Peng49}
V. Kotlyarov,  A. Minakov, Riemann-Hilbert problem to the modified Korteveg-de Vries equation: Long-time dynamics of the steplike initial data, J. Math. Phys. 51 (2010)  093506.
\bibitem{Peng50}
J. Xu,  E. G. Fan,  Y. Chen,  Long-time Asymptotic for the derivative nonlinear Schr\"{o}dinger
equation with step-like initial value, Math. Phys. Anal. Geom. 16 (2013) 253-288.
\bibitem{Peng51}
M. Borghese,  R. Jenkins,  K. T. R. McLaughlin. Long time asymptotic behavior of the focusing nonlinear Schr\"{o}dinger equation. Annales de l'Institut Henri Poincar\'{e}-Analyse non lin\'{e}aire,  35(4) (2018) 887-920.
\bibitem{Peng52}
R. Jenkins,  J. Liu,  P. Perry, et al. Soliton resolution for the derivative nonlinear Schr\"{o}dinger equation. Commun. Math. Phys.,  363(3) (2018) 1003-1049.
\bibitem{Peng53}
Y. L. Yang, E. G. Fan. On the long-time asymptotics of the modified Camassa-Holm equation in space-time solitonic regions. Adv. Math., 402 (2022) 108340.
\bibitem{Tian-jde}
Z. Q. Li,  S. F. Tian,  J. J. Yang, et al. Soliton resolution for the complex short pulse equation with weighted Sobolev initial data in space-time solitonic regions. J. Diff. Equations,  329 (2022) 31-88.
\bibitem{Liunan5}
 G. Biondini,  G. Kova\v{c}i\v{c},  Inverse scattering transform for the focusing nonlinear Schr\"{o}dinger
equation with nonzero boundary conditions. J. Math. Phys. 55 (2014) 031506.
\bibitem{Liunan7}
G. Biondini,  D. Mantzavinos,  Long-time asymptotics for the focusing nonlinear Schr\"{o}dinger
equation with nonzero boundary conditions at infinity and asymptotic stage of modulational
instability, Comm. Pure Appl. Math. 70 (2017)  2300-2365.
\bibitem{Wang}
 D. S. Wang,    B. L. Guo,  X. L. Wang,  Long-time asymptotics of the focusing Kundu-Eckhaus equation with nonzero boundary conditions. J. Diff. Equations, 266(9) (2019)  5209-5253.
\bibitem{Liunan}
B. L. Guo,  N. Liu,  The Gerdjikov-Ivanov-type derivative nonlinear Schr\"{o}dinger equation: Long-time dynamics of nonzero boundary conditions. Math. Meth.  Appl. Sci.,  42(14) (2019) 4839-4861.
\bibitem{Chen-yan}
S. Chen,  Z. Yan,  B. Guo, Long-time asymptotics for the focusing Hirota equation with nonzero boundary conditions at infinity via the Deift-Zhou approach. Math. Phys. Anal.  Geom.,  24(2)(2021) 1-37.
\bibitem{Peng-tian}
W. Q. Peng,  S. F. Tian, Long-time asymptotics in the modified Landau-Lifshitz equation with nonzero boundary conditions. arXiv preprint arXiv:1912.00542, 2019.
\bibitem{Shepelsky}
Y. Rybalko, D. Shepelsky, Asymptotic stage of modulation instability for the nonlocal
nonlinear Schr\"{o}dinger equation. Physica D, 428 (2021) 133060.
\bibitem{Farkas}
H. Farkas,  I. Kra,  Riemann Surfaces, 2nd ed., Grad. Texts in Math., vol. 71, Springer, New York, (1992).
\end{thebibliography}
\end{document}